# The magic of 8 and 24

**Andrei Okounkov**

**Abstract**

While the author is a professional mathematician, he is by no means an expert in the subject area of these notes. The goal of these notes is to share the author's personal excitement about some results of Maryna Viazovska with mathematics enthusiasts of all ages, using maximally accessible, yet precise mathematical language. No attempt has been made to present an overview of the current state field, its history, or to place this narrative in any kind of broader scientific or social context. See the references in Section 5 for both professional surveys and popular science accounts that will certainly give the reader a broader and deeper understanding of the material.

## 1. Spheres keep their distance

### 1.1. Spheres in a $d$-dimensional space

High-dimensional spaces really exist. A photo of a 3-dimensional object taken by our phone may seem to be a 2-dimensional representation of the original, but ... As we capture, process, store, transmit, or display photos, real manipulations are happening with long list of numbers $(x_1, \ldots, x_d)$. Not just inside the phone, but also inside our brain, the image is processed out of many millions of electric potential readouts from the cone and rod cells.

We will call a list of numbers $\boldsymbol{x} = (x_1, \ldots, x_d)$ a vector. Possible values of each $x_i$ may be different in different contexts. It could be just a bit, meaning $x_i$ equals either 0 or 1. It could take values from 0 to 255, as in many popular color specifications. If $x_i$ records a value of the electric potential then, in principle, it is a real number that can take any value, arbitrarily small or large. While these different contexts all influence and enrich each other, our focus in this narrative will be on real vectors. Mathematicians denote real numbers by $\mathbb{R}$ and $d$-tuples of real numbers by $\mathbb{R}^d$. The number $d$ is called the dimension.

To $\mathbb{R}^2$ and $\mathbb{R}^3$, we can attach a familiar geometric image. Via Cartesian coordinates, a point $\boldsymbol{x} = (x_1, x_2) \in \mathbb{R}^2$ corresponds to a point in the plane, whereas $\boldsymbol{x} = (x_1, x_2, x_3) \in \mathbb{R}^3$ corresponds to a point in the our native 3-dimensional space. While $\mathbb{R}^d$ may not be as familiar, it exists and it is important. With diverse uses and applications in mind, mathematicians, scientists, and engineers are all learning to wrap their 3-dimensional heads around the $d$-dimensional spaces.

A key geometric quantity in $\mathbb{R}^3$ is the distance between two points

$$\|\boldsymbol{x} - \boldsymbol{y}\| = \sqrt{\sum_{i=1}^{d}(x_i - y_i)^2}, \tag{1}$$

where $d = 3$. For $d = 2$, this is the distance between two points in the plane. For any $d$, this is the most natural way to define the distance between two points in $\mathbb{R}^d$. It is an important and useful notion in countless contexts, for instance, in statistical analysis.

For example, suppose we measured the values $\boldsymbol{x}' = (x'_1, \ldots, x'_d)$ where we expected to see $\boldsymbol{x} = (x_1, \ldots, x_d)$. Should we attribute the discrepancy to a small unavoidable random noise? Or have we observed something unexpected? The distance $\|\boldsymbol{x}' - \boldsymbol{x}\|$ is the principal measure of how well our measurements fit our predictions.

In this and other situations, it becomes important to separate the points $\boldsymbol{x}'$ whose distance from $\boldsymbol{x}$ is larger than some fixed threshold. One thus defines the ball and the sphere in $\mathbb{R}^d$ with center $\boldsymbol{x}$ and radius $r$, respectively, by

$$B(\boldsymbol{x}, r) = \{\boldsymbol{x}', \text{ such that } \|\boldsymbol{x}' - \boldsymbol{x}\| \leq r\}, \tag{2}$$

$$S(\boldsymbol{x}, r) = \{\boldsymbol{x}', \text{ such that } \|\boldsymbol{x}' - \boldsymbol{x}\| = r\}. \tag{3}$$

We will use this terminology in all dimensions, even though for $d = 2$ this is usually called a disc and circle[1], respectively.

---

[1] And for $d = 1$, (2) is a segment and (3) are its endpoints.

The principal question for us in this narrative is how densely can one pack the spheres of a fixed radius in the $d$-dimensional space. One may compare and contrast a sphere with a cube

$$\text{Cube}(\boldsymbol{x}, r) = \{\boldsymbol{x}', \text{ such that } \max_{i} |x'_i - x_i| \leq r\}, \tag{4}$$

with center $\boldsymbol{x}$ and size $(2r) \times \cdots \times (2r)$. The maximum in (4) is an alternative measure of proximity of two vectors $\boldsymbol{x}$ and $\boldsymbol{y}$, and it is useful in different contexts. Spheres are exceptionally symmetric, preserved by all possible rotations around their center. Compared with spheres, cubes look heavy and boxy. Stacked side to side, cubes fill the whole space, leaving no voids between them. Two spheres can only touch at a point, and there will be voids left no matter how cleverly we try to pack them. However, what is the densest packing that can be achieved?

We will see that different dimensions vary significantly when it comes to sphere packings. In particular, in $\mathbb{R}^8$ and $\mathbb{R}^{24}$ there exist very special arrangements of spheres, denoted $E_8$ and $\Lambda_{24}$. They have been conjectured to be the densest possible in these dimensions.

Recently, this conjecture was proven in an absolutely stunning fashion by Maryna Viazovska in a solo work [54] for $E_8$ and by Viazovska and collaborators Henry Cohn, Abhinav Kumar, Stephen D. Miller, and Danylo Radchenko for $\Lambda_{24}$ in [11]. For these and other phenomenal results, Maryna Viazovska was awarded the Fields Medal, the highest honor in mathematics, in 2022. Our modest goal in these notes is to share our personal excitement about the amazing math that goes into both the statement and the proof of these theorems with the broadest possible audience of mathematics enthusiasts.

**1.2. Sphere packings in $\mathbb{R}^2$**

The problem of sphere packing in 2 dimensions is familiar to anyone who tried to cut circular pieces from a rolled dough while preparing any of the delicious variations on the same universal theme, from vareniki to empanadas [56].

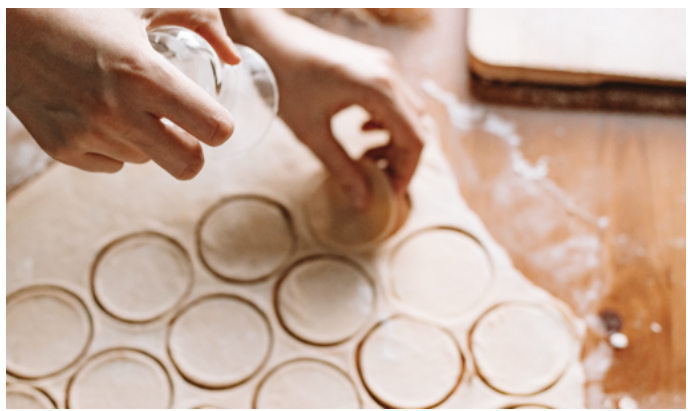

(5)

Naturally, one would like to minimize the fraction of the dough that gets discarded. If the size of the dough is much larger than the radius $r$ of the cutter then this fraction is much more sensitive to the arrangement of circles than to the value of $r$. As a mathematical abstraction, one can consider an infinite piece of dough, and compute the fraction of the dough used (that is, the density of the sphere packing) as a *limit*[2] over larger and larger squares like in Figure (6). We will compute the densities in Figure (6) momentarily.

Note that for the infinite plane, a simple rescaling shows that the packing density does not depend on the radius $r$. This is true for sphere packing in all dimensions. In the analysis, one may leave $r$ as a variable, or set it to any convenient value.

Let us look at Figure (6) more closely:

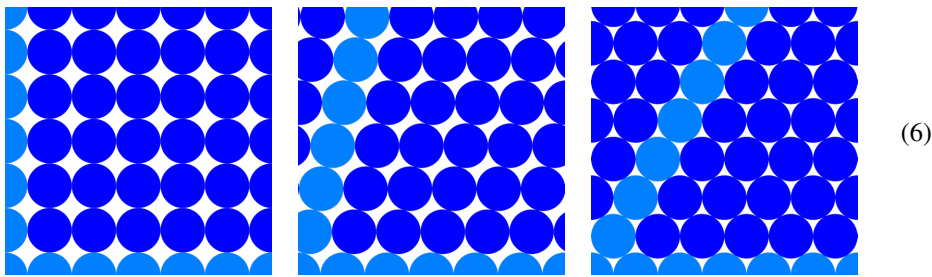

(6)

On the left, we have stacked the circles just like squares. Hence, within each square, it is the inscribed circle that is used, and the rest discarded. Therefore, the packing density is

$$\frac{\text{area of inscribed circle}}{\text{area of a square}} = \frac{\pi}{4} = 0.785\ldots. \qquad (7)$$

If we slant the packing we can improve this. The other two arrangements in Figure (6) are slanted at the angle of $\frac{5\pi}{12} = 75°$ and $\frac{\pi}{3} = 60°$, respectively. Therefore the distance between the horizontal rows of circles has decreased, and namely by a factor of

$$\sin\frac{5\pi}{12} = 0.965\ldots, \quad \sin\frac{\pi}{3} = \frac{\sqrt{3}}{2} = 0.866\ldots.$$

As the horizontal rows get closer, the density increases by the reciprocals of these numbers. At $\frac{\pi}{3}$ this improvement has to stop, because each circle now touches not 4 but 6 other circles, and we cannot slant the figure any further.

Arguably, the hexagonal arrangement on the right in Figure (6), with its 6-fold symmetry, is even more symmetric than the square arrangement on the left. It has a special name in mathematics, namely $A_2$. Here 2 stands for the dimension and the letter $A$ will be discussed a bit later. We have

$$\text{density}(A_2) = \frac{\pi}{2\sqrt{3}} = 0.906\ldots. \qquad (8)$$

This is the densest the spheres can be packed in two dimensions. It is not simple to give a rigorous mathematical proof of this fact, but mathematicians succeeded a long time ago; see [21, 22, 53]. The person cutting the dough in the photograph (5) is evidently aware of this.

---

[2] Readers unfamiliar with limits may probably find their discussion in [43] useful. To avoid worrying about the existence of the limit, it is a good idea to replace the limit by limit superior in this definition.

## 1.3. Contact number in $\mathbb{R}^3$

Let's see how well the life in three dimensions has prepared us for the analysis of the sphere packings in $\mathbb{R}^3$. As a warm-up, one can consider a local version of the packing problem, known as the contact number problem. It can be asked in any dimension and asks for the maximal number $\tau(d)$ of spheres of radius $r$ in $\mathbb{R}^d$ that can be brought in contact with a given sphere of the same radius.

It is quite clear and will be revisited below that $\tau(2) = 6$, realized by the $A_2$ arrangement. In $\mathbb{R}^3$, the problem has a long history, the origin of which legend attributes to the notes taken by David Gregory during his conversations with Isaac Newton in 1694. See [5] for a critical analysis of this legend[3].

Newton and Gregory apparently talked about celestial bodies, in which context it is natural to ask which percentage of the sky on one body, say the Earth, is occupied by image of another body, say the Moon, like in Figure (9):

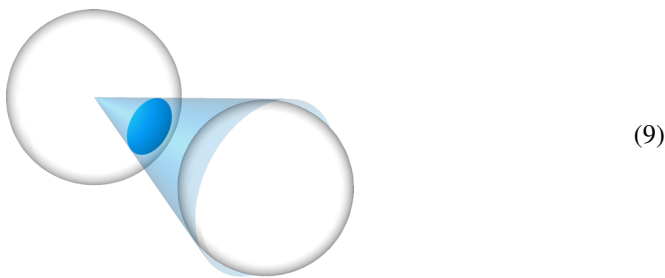

(9)

This percentage depends on the ratio of the distance to the Moon to Moon's radius. Suppose we have two touching spheres of the same radius $r$ like in Figure (10):

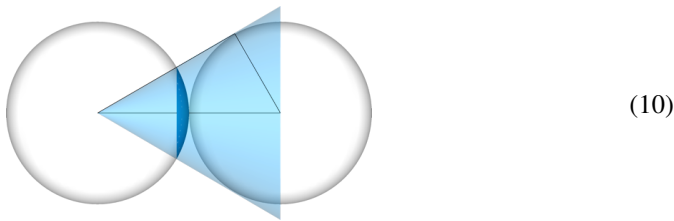

(10)

The right triangle in Figure (10) has hypotenuse $2r$ and short side $r$. Therefore, the opposing angle equals

$$\arcsin \frac{1}{2} = \frac{\pi}{6}, \qquad (11)$$

regardless of the dimension. Note that for $\mathbb{R}^2$ this already suffices to conclude $\tau(2) = 6$.

For $d = 3$, which fraction of the sky is occupied by the spheres in Figure (10) in each other's sky? Consider the sphere

$$S(0, r) = \{(x_1, x_2, x_3) : x_1^2 + x_2^2 + x_3^2 = r^2\} \subset \mathbb{R}^3 \qquad (12)$$

---

[3] It is a problem in mathematics and human life in general that, lacking the time and resources to research every single topic, we mostly just repeat what we have been told. Not being able to break with this tradition, the narrator cannot do better than repeat what he read in [5].

with center at the origin and radius $r$. The points with $x_3 \geq h$, where $h$ is some fixed number between $-r$ and $r$, form what is called a *spherical cap*. The images in the sky in Figures (9) and (10) are spherical caps[4].

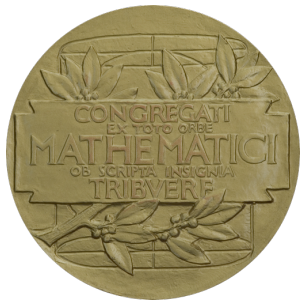

It was known already to Archimedes, and is commemorated as the comparison of the sphere with the cylinder on the back of the Fields Medal, that the area of a spherical cap is *proportional* to its height $r - h$.

Since the cap vanishes for $h = r$ and is the whole sphere for $h = -r$, we conclude

$$\frac{\text{area of the cap}}{\text{area of the sphere}} = \frac{r-h}{2r} = \frac{1 - h/r}{2}. \qquad (13)$$

For the cap in Figure (10),

$$\frac{h}{r} = \cos\frac{\pi}{6} = \frac{\sqrt{3}}{2} \quad \Rightarrow \quad \frac{2r}{r-h} = \frac{4}{2-\sqrt{3}} = 14.92\ldots \qquad (14)$$

Since each cap occupies more that $1/15$ of the surface area, 15 caps cannot fit without overlap, so the contact number in $\mathbb{R}^3$ is at most 14. It is easy to see that it is at least 12. Can it be equal 13? In the legend, Gregory thought yes, while Newton thought no.

There is an objective difficulty here, and it has to do with the fact that there are many different possible configurations of 12 spheres. One of them, realized for the densest packing, can be seen on the left in Figure (27) below. But another possibility is to put the spheres in the 12 vertices of a regular icosahedron, like in Figure (15), which also shows the plot of the spherical caps.

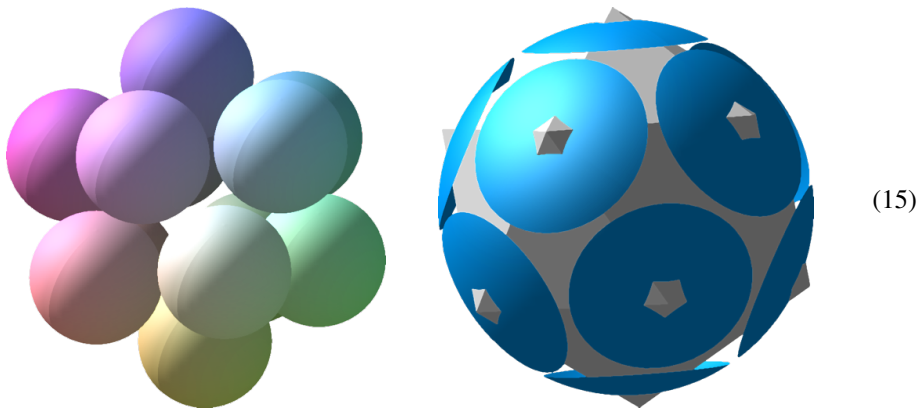

(15)

From the spherical caps, we see that the 12 spheres are not touching each other, hence can be moved around without losing the contact with the central sphere[5]. Perhaps we can make room for the 13th sphere? The problem of fitting caps into a sphere is a close spherical relative of the sphere packing problem[6] and belongs to a broad class of problems known as spherical codes and spherical designs [2, 17].

---

4     Of course, one gets the same geometric shape if one takes $x_1 \geq h$ or $x_2 \geq h$ instead of $x_3 \geq h$. Since caps are normally worn on the top of the head, we made the conventional choice of the vertical $x_3$-axis.

5     It is a fun fact to prove, see the appendix to Chapter 1 in [15], that an *arbitrary* permutation of the 12 spheres may be achieved by rolling them around the central sphere.

Popular descriptions of the contact number problem often contain a suggestion for the reader to imagine a billiard ball hanging in mid-air, and 12 further billiard balls rolling around it. I envy those readers who have enough spatial intuition to imagine something like this. Even in our native $\mathbb{R}^3$, our geometric intuition often asks for help. Help may come from building models or from doing computations.

Geometers of all times have liked building models, using whatever materials the technology of the time made available. They would be surely thrilled to see the computer models that we can build today. There is a wonderful animated popular account of the contact number problem at Mathematical Etudes website. I am sure many readers will find it fascinating.

But the destination of our story being sphere packings in $\mathbb{R}^8$, it may be safe to expect the computations to overtake models in such high dimension. And indeed, at the heart of Viazovska's proof in [54] lies a brilliant inspired computation. It puts a big exclamation mark in a certain long line of argument. This line of argument was first born in the work of Philippe Delsarte in the discrete setting of coding theory [16] and was later adapted to spherical codes to compute

$$\tau(4) = 24, \quad \tau(8) = 240, \quad \tau(24) = 196560, \quad (16)$$

see [38, 40, 42]. It may be also used to show that $\tau(3) = 12$, see [1, 40], but many other proofs of this fact were found earlier [5, 47]. The Newton character from the legend was right.

Delsarte-type bounds, also known as linear programming bounds, were put to work in the sphere packings situation by H. Cohn and N. Elkies in [10]. We will talk about them in Section 3. They require a certain magic function to complete the proof. It is this elusive magic function that was discovered by Viazovska is her astonishing work [54].

---

[6] Strictly speaking, taking into account the spherical shape of the Earth, the person in Figure (5) may be solving the spherical cap packing problem. In contrast to the sphere packing problem in $\mathbb{R}^d$, the radius of the spherical cap, or equivalently its angular size, cannot be scaled away and remains an important parameter in the problem. In the limit of very small caps, the problem reduces to sphere packing in the flat space $\mathbb{R}^d$. Of course, a person making vareniki is not taking the radius of the Earth into account!

## 1.4. The densest packings in $\mathbb{R}^3$

Let

$$v_1, \ldots, v_d \in \mathbb{R}^d$$

be a basis of $\mathbb{R}^d$, equivalently a set of linearly independent vectors[7]. By definition, the *lattice* $\Lambda$ spanned by the vectors $v_1, \ldots, v_d$ is formed by all vectors

$$\Lambda = \mathbb{Z}v_1 + \mathbb{Z}v_2 + \cdots + \mathbb{Z}v_d \subset \mathbb{R}^d \tag{17}$$

that can be obtained from the $v_i$'s using addition and subtraction. The linear space $\mathbb{R}^d$ is a group under addition and lattices are special kinds of subgroups in it.

A sphere packing is called a *lattice packing* if the centers of the spheres form a lattice. For instance, the packings in Figure (6) are lattice packings. There, we can take $v_1 = (2r, 0)$ is all three cases, while

$$v_2 = (2r\cos\phi, 2r\sin\phi), \quad \text{where} \quad \phi = \tfrac{\pi}{2}, \tfrac{5\pi}{12}, \tfrac{\pi}{3},$$

depending on the slant angle $\phi$.

Both the hexagonal packing and the corresponding hexagonal lattice are denoted by the symbol $A_2$. There is a cool way to realize this lattice inside $\mathbb{R}^3$ as the set

$$A_2 = \{(x_1, x_2, x_3), \sum x_i = 0\} \subset \mathbb{Z}^3 \tag{18}$$

of integer points with sum zero; see Figure (19):

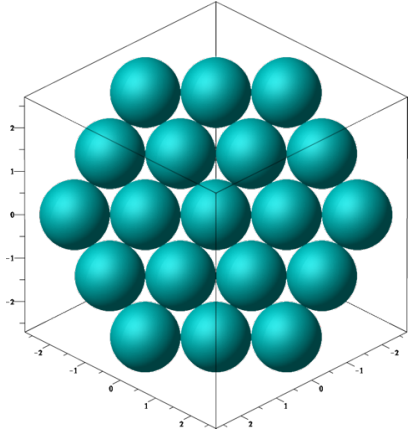

(19)

In (19), the sphere at $(0, 0, 0)$ is surrounded by 6 spheres with centers at all possible permutations of $(1, -1, 0)$. These are at distance $\sqrt{2}$ from the origin, and hence $r = \tfrac{1}{\sqrt{2}}$.

In a second we will need to talk about holes in the $A_2$ packing, shown in Figure (21). These come in two different flavors according to the sign in

$$\text{hole center} = \pm(\tfrac{2}{3}, -\tfrac{1}{3}, -\tfrac{1}{3}) + \text{integer vector}. \tag{20}$$

---

[7] Some readers may find the explanation of these notions given in [44] useful.

The two kinds of holes are color-coded in Figure (21). They are permuted by symmetries of $A_2$.

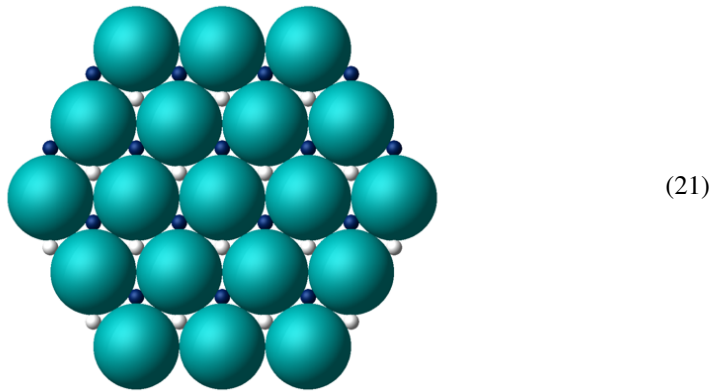

(21)

The densest sphere packing in $\mathbb{R}^3$ may be constructed by adding new layers of spheres (19) as in Figure (22). Each new layer is a copy of (19) shifted so that the new spheres fit over the holes of the previous layer

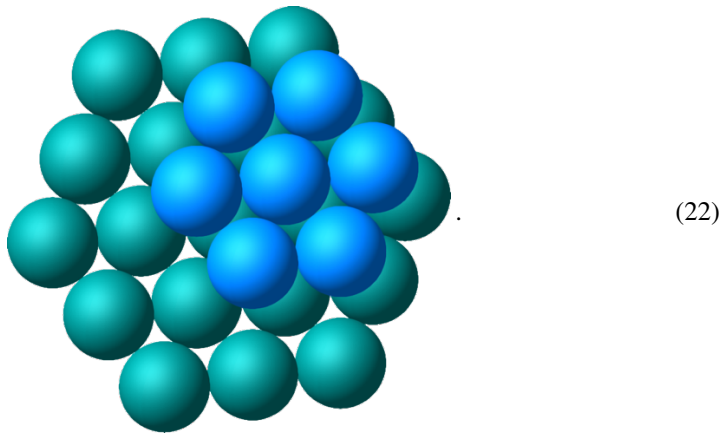

. (22)

Since at every step we have 2 possible choice of the holes in (20), this gives $2^\infty$ different choice of packings with the same density! However, if we want it to be a lattice packing then there is only one choice up to an overall rotation or reflection. We can take

$$\begin{aligned} v_1 &= (1,-1,\ 0), \\ v_2 &= (0,\ 1,-1), \\ v_3 &= (0,\ 1,\ 1), \end{aligned} \qquad (23)$$

and these generate the lattice

$$D_3 = \{(x_1, x_2, x_3), \sum x_i \text{ is even}\} \subset \mathbb{Z}^3. \qquad (24)$$

In general, one defines

$$A_d = \{(x_1, \ldots, x_{d+1}), \sum x_i = 0\} \qquad \subset \mathbb{Z}^{d+1}, \qquad (25)$$

$$D_d = \{(x_1, \ldots, x_d),\ \sum x_i \text{ is even}\} \quad \subset \mathbb{Z}^d. \qquad (26)$$

For $d > 3$ these define different lattices and different sphere packings, but it is mathematical fact that $A_3$ is the same as $D_3$. Check this! Henry Cohn suggests the following exercise for when the reader visits the grocery store next. Find some fruit stacked as $A_3$ and some stacked as $D_3$. Then rotate your head until you are convinced that they are the same packing!

The proof of the fact that $D_3$ is the densest sphere packing in $\mathbb{R}^3$ is a monumental achievement of T. Hales and an inspiring story of computers helping humans to finish very complex proofs. See [25, 27, 28, 37] for more about this. It is not known, but conjectured, that $D_4$ and $D_5$ are the densest sphere packings in $\mathbb{R}^4$ and $\mathbb{R}^5$, respectively.

It should be stressed emphatically that the optimality discussed in these notes concerns optimality among all sphere packings, not just lattice packings. Within the class of lattice packings, the optimality of $D_3$ was shown by Gauss in 1831 [23], while the optimality of the $D_4$ and $D_5$ lattices was proven by Korkine and Zolotareff [33, 34] in 1870s.

Following our discussion of the contact number $\tau(3)$, it is fun to examine the arrangement of neighbors in the $A_3 = D_3$ packing. The sphere at $(0, 0, 0)$ has 12 neighbors with centers at the vectors $(\pm 1, \pm 1, 0)$ and their permutations. These are the 12 vectors $x$ of length $\|x\| = \sqrt{2}$ in the lattice $A_3 = D_3$. The corresponding spherical caps can be seen in Figure (27), together with 24 caps for the spheres with centers at points $\|x\|^2 = 6$ and 48 caps for the spheres at the distance $\|x\|^2 = 14$

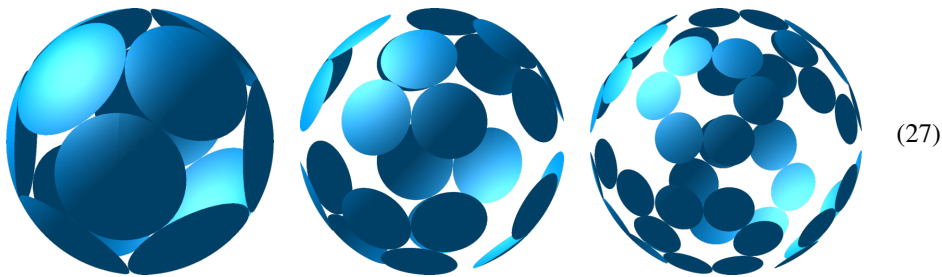

(27)

In this fashion, one can obtain very interesting collections of points on spheres from dense lattice packings in any dimension.

## 2. Beyond the 3-space
### 2.1. 4, 5, 6, 7, 8, …

The narrator of these notes is a complete novice in the field of sphere packing trying to share his first impressions of the striking beauty of the field with other mathematics enthusiasts. Among the mathematicians of older generations, I imagine I am not alone feeling like a schoolboy again, exploring spellbound the treasures described, in particular, in the treatise [15] by John Conway, Neil Sloane, and collaborators. The story starts deceptively simple but quickly leads to the highest heights and deepest depths of mathematics.

To continue the parallel with one's student years, each dimension $d$ in the sphere packing problem feels like a new year of math classes. While it builds on and connects with the material form the previous years, many new phenomena and ideas appear each time.

In even further parallel to how mathematics courses change as we go trough high-school, college, graduate school, and so on, hopefully never stopping learning, the sphere packing problems seem to come in certain groups of dimensions. In dimensions up to 8, the densest packing are known or conjectured to be the lattice packings

$$\underline{A_1}, \underline{A_2}, \underline{A_3}, D_4, D_5, E_6, E_7, \underline{E_8}. \tag{28}$$

Here the known cases are underlined, the 8-dimensional case being Viazovska's breakthrough. This certainly feels like a story about exceptional Lie groups, which ends in dimension 8 with the largest exceptional Lie group $E_8$.

Next come dimensions 9 through 24. Looking at the iconic picture in (29), reproduced here with permission from [49], we see that these dimensions start out as valley leading to an ascent to the sharp peak of the Leech lattice $\Lambda_{24}$. The Leech lattice is now proven to give the densest sphere packing in dimension 24 by Viazovska and collaborators [11].

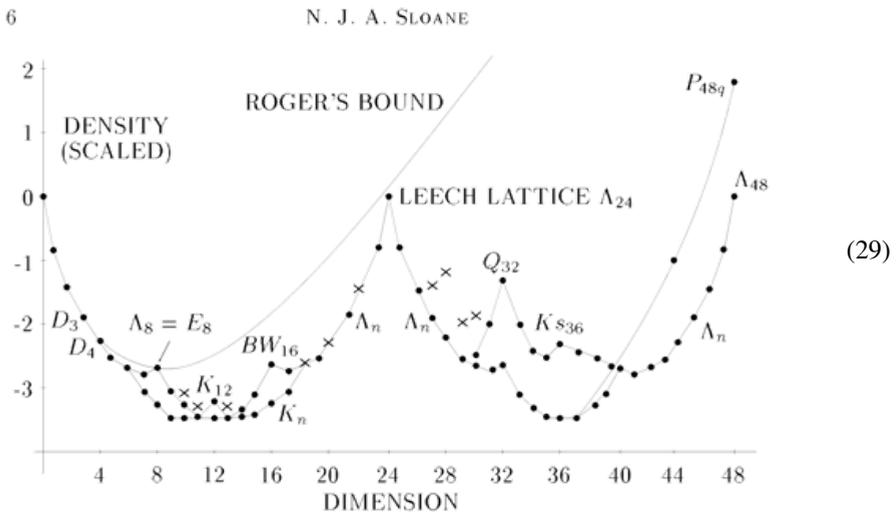

(29)

The Leech lattice, with its deep connections to the exceptional, or sporadic, finite simple groups including the Monster group of Bernd Fischer and Robert Griess, is the defining feature of the 9 to 24 valley. Again, the Monster being the largest exceptional group, the storyline has to change after 24.

What is next? We hope the reader's curiosity will lead her or him to explore, guided by [6, 15]. See also the tables [9, 41] of the densest sphere packings currently known in different dimensions.

### 2.2. Fluid diamond in $d = 9$

Here is one among the countless marvels of high-dimensional sphere packing. Recall the lattice $D_d$ from (25) formed by integer vectors with even coordinate sum. The nearest neighbors in $D_d$ are $\sqrt{2}$ away, so we can pack spheres of radius $r = \frac{1}{\sqrt{2}}$ using points of $D_d$ as centers.

Consider the vector
$$\boldsymbol{\gamma} = (\tfrac{1}{2}, \tfrac{1}{2}, \tfrac{1}{2}, \ldots, \tfrac{1}{2}) \,. \tag{30}$$

What is the distance between $x$ and the nearest point $v \in D_d$? Since $v$ has integral coordinates, we have
$$\|\boldsymbol{\gamma} - \boldsymbol{v}\|^2 \geq \underbrace{\frac{1}{2^2} + \cdots + \frac{1}{2^2}}_{d \text{ times}} = \frac{d}{4} \,.$$

Therefore, if $d = 8$, we can fit *two* copies of $D_8$ into $\mathbb{R}^8$ with a shift by $\boldsymbol{\gamma}$. The resulting lattice is nothing else than the magic $E_8$ lattice
$$E_8 = D_8 \cup (D_8 + \boldsymbol{\gamma}) \,, \tag{31}$$

about which we will talk more in Section 2.3 below.

If $d = 9$, we can take the vector
$$\boldsymbol{\gamma}_{i,t} = \boldsymbol{\gamma} + t\boldsymbol{e}_i \,, \quad i = 1, \ldots, 9 \,, \tag{32}$$

where $t \in \mathbb{R}$ is an arbitrary number and
$$\boldsymbol{e}_i = (0 \ldots, 0, \underset{i}{1}, 0, \ldots, 0) \tag{33}$$

is the $i$th coordinate vector. By the same argument as before
$$\|\boldsymbol{\gamma}_t - \boldsymbol{v}\|^2 \geq 2 \,, \quad \text{for} \quad \boldsymbol{v} \in D_9 \,,$$

so we can pack the spheres of radius $r = \frac{1}{\sqrt{2}}$ using points of
$$\text{fluid diamond packing} = D_9 \cup \left( D_9 + \boldsymbol{\gamma}_{i,t} \right) \tag{34}$$

as centers. Note that since both $t$ and $i$ are arbitrary, half of the spheres in (34) can be shifted arbitrarily in one of the coordinate direction without running into the other half of the spheres — a rather fluid packing! And yet, its density matches, for any $t$, the highest known density in dimension 9. It is not so easy to imagine this possible based on our low-dimensional geometric intuition.

### 2.3. Stars align in $E_8$

The exceptionally dense and symmetric $E_8$ lattice packing which we met in (31) certainly merits a much longer discussion. One can start this discussion from many different angles, emphasizing different areas of mathematics where the $E_8$ lattice naturally appears.

#### 2.3.1. Roots

The lattices $A_d$ and $D_d$ from (25) have the property that $\|v\|^2$ is an even integer for any $v \in D_d$. Such lattices are called *even*. How can we tell if a lattice $\Lambda$ as in (17) is even? Using the concept of the *inner product*, recalled in Appendix A, it suffices to check that $(v_i, v_j) \in \mathbb{Z}$ and $(v_i, v_i) \in 2\mathbb{Z}$ for any basis of $\Lambda$. Since

$$(\gamma, \gamma) = 2 \quad \text{and} \quad (\gamma, v) \in \mathbb{Z},$$

for any $v \in D_8$ and $\gamma$ as in (30), we see that $E_8$ is an even lattice.

Given an even lattice $\Lambda$, vectors $\alpha \in \Lambda$ of the minimal nonzero norm $\|\alpha\|^2 = 2$ are call *roots*[8]. These are the centers of the spheres touching the central sphere. For example, the vectors

$$\alpha = \pm e_i \pm e_j \in D_d, \quad i \neq j, \tag{35}$$

are roots. For $E_8$, we also have the root $\gamma$ as well as

$$\alpha = (\pm\tfrac{1}{2}, \pm\tfrac{1}{2}, \pm\tfrac{1}{2}, \ldots, \pm\tfrac{1}{2}), \quad \text{such that the sum is even,} \tag{36}$$

of which there are $\tfrac{1}{2}2^8 = 128$ many. We invite the reader to check there are no other roots for $E_8$ and verify that the number of roots equals $\tau(8) = 240$. Thus the roots of $E_8$ give the solution of the contact number problem in $d = 8$, and in fact this solution is unique, very much unlike the $d = 3$ case discussed in Section 1.3.

#### 2.3.2. Reflections

Every root $\alpha \in \Lambda$ in an even lattice $\Lambda$ generates a special symmetry of the lattice $\Lambda$, namely the orthogonal reflection $r_\alpha$ in the hyperplane orthogonal to $\alpha$. It sends $\alpha$ to $-\alpha$ and fixes all vectors $v_\perp$ that are orthogonal to $\alpha$.

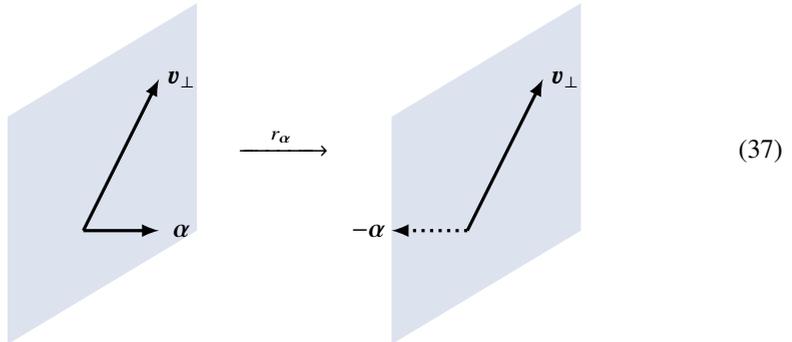

(37)

---

[8] For lattices in which norm squared takes both even and odd integer values, vectors of norm 1 should also count as roots. These have an important role to play in Lie theory and many other branches of mathematics, but not in our narrative.

Explicitly, it is given by the formula

$$r_\alpha(v) = v - (v, \alpha)\,\alpha,\tag{38}$$

which manifestly preserves the lattice $\Lambda$. Indeed, (101) shows the inner product takes integer values in an even lattice. For more on (38), see Section A.3.

For example, the roots $\alpha$ of the $A_2$ lattice (18) are the permutations of the vector $(1, -1, 0)$. The reflection $r_{(1,-1,0)}$ swaps the first two coordinates. Similarly, for $A_d \subset \mathbb{R}^{d+1}$, each reflection $r_\alpha$ swaps two coordinates of $\mathbb{R}^{d+1}$.

Orthogonal symmetries of a lattice $\Lambda$ always form a finite group; see the brief introduction to this concept in Appendix B. In particular, its subgroup generated by the reflections $r_\alpha$ is a finite group $W$ generated by reflections. Such groups have been fully classified and studied in great detail due to their crucial importance in Lie theory, singularity theory, and many other branches of mathematics. As a corollary of this classification, we know that all even lattices spanned by roots are orthogonal direct sums of lattices of the form $A_d$, $D_d$, or $E_6, E_7, E_8$.

### 2.3.3. ADE classification

How does this classification work? In an even lattice $\Lambda$ spanned by roots one can always choose the basis of roots so that

$$(\alpha_i, \alpha_j) = 0 \text{ or } -1, \quad i \neq j.\tag{39}$$

For example, for $E_8$, we can take

$$\alpha_i = e_i - e_{i+1}, \quad i = 1, \ldots, 6,\tag{40}$$

together with

$$\alpha_7 = e_6 + e_7, \quad \alpha_8 = -\gamma.\tag{41}$$

We see that $(\alpha_i, \alpha_j) = 0$ for most pairs $i, j$, and that $(\alpha_i, \alpha_j) = -1$ precisely for the pairs connected by an edge in the following graph:

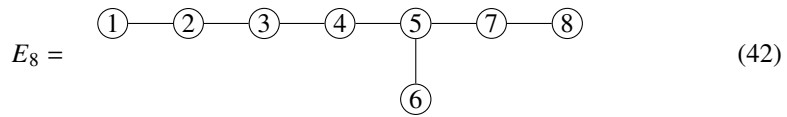

$$E_8 = \tag{42}$$

The graph (42) is a very convenient graphical way to represent the Gram matrix

$$B_\Lambda = ((\alpha_i, \alpha_j)).$$

It is called the Dynkin diagram or the Coxeter diagram.

The ADE classification is really the classification of all possible diagrams like (42) for which the corresponding Gram matrix $B_\Lambda$ is positive definite (a concept which will be explained and used in Section 3.1 below). This is not as difficult as it sounds, and revolves around the fact any subgraph of a positive definite diagram is a positive definite diagram.

For instance, if we erase nodes $1, 2, 3, 8, 7, 6, 5$ from (42) in that order, we get diagrams for lattices $E_7, E_6, D_5, D_4, D_3 = A_3, A_2, A_1$, that are known or conjectured to give the densest packings in the corresponding dimensions.

In the opposite direction, if we try to invent the lattice $E_9$ with the diagram

$$\text{``}E_9\text{''} = \begin{array}{c} \text{(0)}-\text{(1)}-\text{(2)}-\text{(3)}-\text{(4)}-\text{(5)}-\text{(7)}-\text{(8)} \\ | \\ \text{(6)} \end{array} \tag{43}$$

we see that this doesn't work because the determinant $\det B_{E_9}$ of the corresponding Gram matrix vanishes[9].

### 2.3.4. Discriminant

In general, the determinant or the discriminant of an even lattice $\Lambda$

$$\Delta_\Lambda = \det B_\Lambda \tag{44}$$

is a very important quantity that enters the following formula for the density of the corresponding sphere packing.

The $\Lambda$-translates of the parallelepiped

$$\Pi_\Lambda = \left\{ x = \sum x_i \alpha_i, \ \max |x_i| \le \tfrac{1}{2} \right\} \tag{45}$$

tile the whole space, each tile containing exactly one lattice point as its center. This just the story about the cube (4) in different coordinates. We have

$$\operatorname{Vol} \Pi_\Lambda = \sqrt{\Delta_\Lambda} \,. \tag{46}$$

Therefore, if $\Lambda$ is even and the roots are the vectors of minimal length then

$$\text{density of the sphere packing} = \frac{\operatorname{Vol} B(0, \frac{1}{\sqrt{2}})}{\sqrt{\Delta_\Lambda}} \,. \tag{47}$$

The smaller the discriminant the larger the density. Since the discriminant is an integer, 1 is the smallest it can be, and in fact

$$\Delta_{E_8} = 1 \,. \tag{48}$$

Lattices with $\Delta_\Lambda = 1$ are called *unimodular*. Even unimodular lattices exist only in dimensions that are multiple of 8, and $E_8$ is the unique even unimodular lattice in $\mathbb{R}^8$.

In dimension 24 there exist 24 even unimodular lattices, and the superamazing Leech lattice is distinguished among them by having no roots! In other words, one can pack spheres radius $r = 1$ with centers in Leech lattice instead of $r = \frac{1}{\sqrt{2}}$. While a meaningful discussion of the Leech lattice transcends the introductory nature of these notes, we hope that the reader's curiosity will be satisfied by the accounts in [15, 18, 52].

---

[9] In fact, the diagram (43) is the Dynkin diagram of the infinite affine reflection group of type $\widehat{E}_8$ in $\mathbb{R}^8$.

### 2.3.5. Codes

Recall how at the very beginning, in Section 1.1 we talked about the possible values of the entries $x_i$ in a vector $\boldsymbol{x} = (x_1, \ldots, x_d)$. While everywhere else in this narrative we consider the case of real entries $x_i$, let's turn our attention to the case $x_i \in \{0, 1\}$ for a brief moment. In other words, let's talk about binary vectors.

The natural distance between binary vectors is the Hamming distance

$$\|\boldsymbol{x} - \boldsymbol{x}'\|_{\text{Hamming}} = \sum |x_i - x_i'|. \tag{49}$$

It measures the number of entries in which $\boldsymbol{x}$ and $\boldsymbol{x}'$ differ, and it is very natural for error correction and other applications. If $\boldsymbol{x}$ and $\boldsymbol{x}'$ represent the input and output of a transmission through a binary communication channels, then (49) is the number of errors that occurred during the transmission. If we can pack nonintersecting Hamming balls of radius $r$ in $\{0, 1\}^d$ then the centers $\mathscr{C} \subset \{0, 1\}^d$ of these balls give binary code words of length $d$ that corrects up to $r$ errors. A related concept is the minimal distance $\delta$ between the code words from $\mathscr{C}$. Evidently $\delta > 2r$.

Given a code $\mathscr{C}$, we define $\widehat{\mathscr{C}} \subset \mathbb{Z}^d$ as the set of integer vectors that have the same parity as some code word from $\mathscr{C}$. Clearly, if $\boldsymbol{v} \neq \boldsymbol{v}'$ are two distinct points of $\widehat{\mathscr{C}}$ then

$$\|\boldsymbol{v} - \boldsymbol{v}'\| \geq \min(2, \sqrt{\delta}),$$

and hence we can pack sphere of half that radius with centers at $\widehat{\mathscr{C}}$.

For $d = 8$, there exists a remarkable code $\mathscr{C}$ with $\delta = 4$. It is obtained by adding the parity bit to Hamming $(7, 4)$-code, see [1, 15]. The corresponding packing $\widehat{\mathscr{C}}$ is isomorphic to the $E_8$ packing. Similarly, the Leech lattice can be obtained from the Golay code.

### 2.3.6. The Coxeter plane

There is the following cool way to visualize roots for any finite reflection group (requires familiarity with eigenvalues and also with complex numbers, see Section A.5). The material in this section may feel a bit advanced and it could be a good idea to come back to it after reading the material in the Appendix.

Recall the basis $\boldsymbol{\alpha}_i$ from Section 2.3.3 and consider the corresponding reflections $r_{\boldsymbol{\alpha}_i}$. Consider the product $\mathsf{C}$ of all these reflections taken in some order. The reflections do not commute, so $\mathsf{C}$ depends on the order. Remarkably, however, the conjugacy class of $\mathsf{C}$ is independent of the order. A a particularly nice choice is

$$\text{Coxeter element } \mathsf{C} = \underbrace{r_{\alpha_1} r_{\alpha_3} r_{\alpha_5} r_{\alpha_8}}_{\text{commute}} \underbrace{r_{\alpha_2} r_{\alpha_4} r_{\alpha_6} r_{\alpha_7}}_{\text{commute}}, \tag{50}$$

which presents $\mathsf{C}$ as a product of two involutions, that is, two elements that each square to $\mathbb{1}$, where $\mathbb{1}$ is the identity matrix.

The order and the eigenvalues of a Coxeter element can be computed abstractly. For $E_8$, we have $\mathsf{C}^{30} = \mathbb{1}$ and the eigenvalues of $\mathsf{C}$ are exactly the primitive roots of unity of order 30 or, equivalently, the roots of the cyclotomic polynomial

$$z^8 + z^7 - z^5 - z^4 - z^3 + z + 1 = 0. \tag{51}$$

We can take any one of them and project the roots onto corresponding eigenspace $\mathbb{C} \subset \mathbb{C}^8$, called the Coxeter plane. The resulting collection of points are the centers of the circles in Figure (52). The radii in that figure have no exact mathematical meaning and are simply adjusted to resemble a sphere packing. The colors will be explained below.

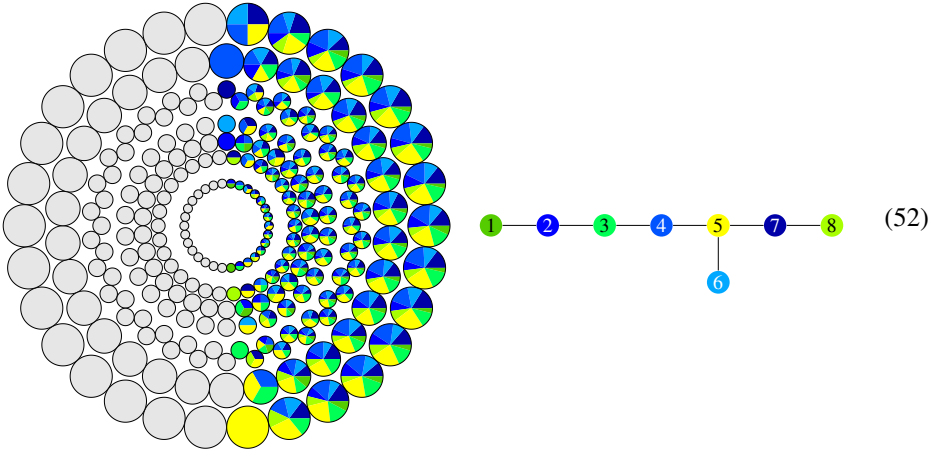

(52)

Consider the involutions

$$\mathsf{R}_{\text{green}} = \underbrace{r_{\alpha_1} r_{\alpha_3} r_{\alpha_5} r_{\alpha_8}}_{\text{commute}}, \quad \mathsf{R}_{\text{blue}} = \underbrace{r_{\alpha_2} r_{\alpha_4} r_{\alpha_6} r_{\alpha_7}}_{\text{commute}}, \qquad (53)$$

where the colors refer to the color-coding in the Dynkin diagram in Figure (52). By construction, the involutions (53) satisfy

$$\mathsf{R}_{\text{green}} \mathsf{C} = \mathsf{R}_{\text{blue}} = \mathsf{C}^{-1} \mathsf{R}_{\text{green}}, \qquad (54)$$

and hence generate, together with $\mathsf{C}$, the symmetry group of a regular 30-gon in the Coxeter plane. We have

$$\mathsf{R}_{\text{green}} \alpha_i = -\alpha_i, \quad i \in \{1, 3, 5, 8\}. \qquad (55)$$

Therefore, all green vertices land on a line in the Coxeter plane — the line perpendicular to the line fixed by $\mathsf{R}_{\text{green}}$. An identical argument works for $\mathsf{R}_{\text{blue}}$.

We can use Figure (52) to illustrate the following important concepts related to roots. First, the roots can be partitioned into positive and negative by a generic hyperplane in $\mathbb{R}^8$, which we can take to be the preimage of a line in the Coxeter plane. In Figure (52), the negative roots are in gray, while the positive roots are colored.

Second, one can choose the roots $\alpha_i$ from Section 2.3.3 as the *simple positive roots*. These are the positive roots that cannot be written nontrivially as a sum of positive roots. They are monochromatic in Figure (52) where the colors correspond to the coloring of the Dynkin diagram as before.

Third, all positive roots are nonnegative integer linear combinations of simple roots. The proportions in which the simple roots combine to produce a given positive root are plotted as pie charts in Figure (52). In particular, the dichromatic roots in Figure (52) correspond to the roots

$$\alpha_i + \alpha_j = r_{\alpha_i}(\alpha_j) = r_{\alpha_j}(\alpha_i) \quad \text{when} \quad (\alpha_i, \alpha_j) = -1, \qquad (56)$$

which exist for any pair of neighbors in the Dynkin diagram.

See Appendix F for more on connections between $E_8$ and regular polygons.

### 2.4. Very large dimensions

Our discussion of sphere packings in arbitrarily large dimensions will be very brief due to both objective lack of information about them and limits of the present narrative.

Let us call a sphere packing in $\mathbb{R}^d$ saturated if no additional sphere of the same radius $r$ can be inserted into it. Remarkably, the density of a a saturated packing is at least $2^{-d}$. We invite the reader to pause for a second and try to prove this. Maryna Viazovska says this is one of her favorite entry-level problems about sphere packings.

One way to prove this is to note that, for a saturated packing, balls of twice the radius with the same centers have to cover the whole $\mathbb{R}^d$. Otherwise, there would a point where we can insert another sphere of radius $r$. From

$$\text{Vol } B(0,r) = 2^{-d} \text{ Vol } B(0, 2r),$$

we get the sought lower bound for the density of a saturated packing.

As simple as this sounds, this bound is remarkable. As we review in Appendix E, the volume of $B(0,r)$ decays superexponentially with dimension $d$ for any $r$. Hence a packing achieving a $2^{-d}$ density must have superexponentially many spheres in any cube $[0, L]^d \subset \mathbb{R}^d$ as $d \to \infty$. Also, the best known improvements to the $2^{-d}$ lower bound are only basically linear in $d$.

As to the upper bounds on density, the one by Kabatiansky and Levenshtein has been holding the world record at $2^{-0.5990...d}$ since 1978, although the methods described below allowed Cohn and Zhao [13] to achieve a constant factor improvement.

## 3. Upper bounds on packing density
### 3.1. Positive definite forms and functions
#### 3.1.1.

Let's start with with simplest possible inequality. For any real number $x$, $x^2 \geq 0$. As trivial as this sounds, this proves, for instance, that

$$x_1^2 - 2x_1x_2 + 2x_2^2 - 2x_2x_3 + x_3^2 = (x_1 - x_2)^2 + (x_2 - x_3)^2 \geq 0. \tag{57}$$

#### 3.1.2.

An expression of the form

$$B(\boldsymbol{x}) = \sum_{i,j=1}^{n} b_{ij}x_ix_j, \tag{58}$$

where $b_{ij} \in \mathbb{R}$ are coefficients, is called a *quadratic form* in the variables $\boldsymbol{x} = (x_1, \ldots, x_n)$. In the sum (58), we may and will assume that $b_{ij} = b_{ji}$. The symmetric array of numbers $(b_{ij})$ is called the matrix of the quadratic form (58).

A quadratic form is called *positive semidefinite* if it takes only nonnegative values, like the one in (57). One writes $B \geq 0$. Forms that take positive values for nonzero arguments are called *positive definite*. For instance, (57) is positive semidefinite but not definite, since it vanishes for $x = (1, 1, 1)$.

### 3.1.3.

If $B_1, B_2 \geq 0$ then
$$c_1 B_1 + c_2 B_2 \geq 0 \tag{59}$$
for all coefficients $c_1, c_2 \geq 0$. Mathematicians say that the set of positive semidefinite forms is a convex cone.

For example, for $n = 2$, we have
$$b_{11}x_1^2 + 2b_{12}x_1x_2 + b_{22}x_2^2 \geq 0 \quad \Leftrightarrow \quad \begin{array}{c} b_{11} \geq 0, \ b_{22} \geq 0 \\ b_{12}^2 \leq b_{11}b_{22} \end{array}. \tag{60}$$

In the 3-space with coordinates $(b_{11}, b_{22}, b_{12})$, the set of the positive semidefinite forms is the familiar cone plotted in Figure (61) with its vertex at the origin $(b_{11}, b_{22}, b_{12}) = (0, 0, 0)$.

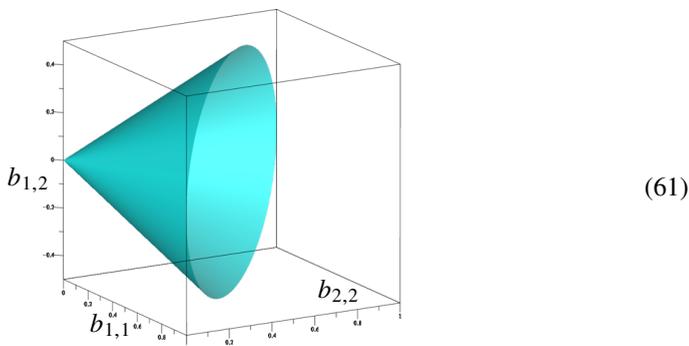

(61)

The interior of this cone corresponds to positive definite forms.

### 3.1.4.

Here is an example of an interesting positive semidefinite form. Fix some angles $\phi_1, \ldots, \phi_n$ and let
$$\varepsilon_i = (\cos \phi_i, \sin \phi_i) \in \mathbb{R}^2, \quad i = 1, \ldots, n,$$
be unit vectors in $\mathbb{R}^2$ having the angle $\phi_i$ with the horizontal axis. Let us add them with coefficients $x_1, \ldots, x_n$; see Figure (62).

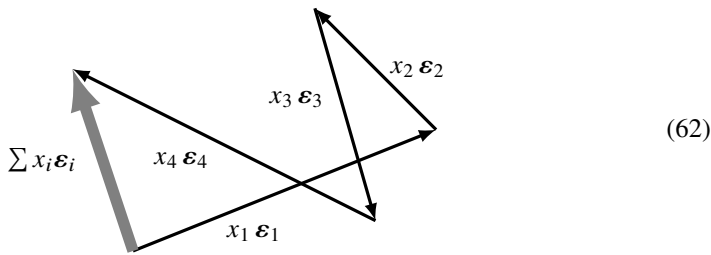

(62)

We define the form $B_{\cos}(\boldsymbol{x})$ as the length squared of this sum. Using inner products, see Appendix A, we compute

$$B_{\cos}(\boldsymbol{x}) = \left\|\sum x_i \boldsymbol{\varepsilon}_i\right\|^2 = \left(\sum x_i \boldsymbol{\varepsilon}_i, \sum x_j \boldsymbol{\varepsilon}_l\right) =$$
$$= \sum_{ij} (\boldsymbol{\varepsilon}_i, \boldsymbol{\varepsilon}_l)\, x_i x_j = \sum_{ij} \cos(\phi_i - \phi_j)\, x_i x_j \geq 0. \quad (63)$$

This is not obviously nonnegative based on coefficients, but we know it must be nonnegative as a length squared.

### 3.1.5.

In general, let $f(t)$ be an even function; that is, let $f(t)$ satisfy

$$f(t) = f(-t).$$

We say that $f$ is *positive definite* if the quadratic form

$$B_f(\boldsymbol{x}) = \sum_{i,j=1}^{n} f(t_i - t_j)\, x_i x_j \geq 0 \quad (64)$$

is positive semidefinite for *any* choice of $t_1, \ldots, t_n$. It follows from (63) that $f(t) = \cos t$ is positive definite. Similarly, $f(t) = \cos \omega t$ is positive definite for any frequency $\omega$.

To get a better feeling for positive-definite functions, the reader may want to deduce from (60) that $f(0) > 0$ for any positive-definite function $f(t)$ that is not identically zero.

### 3.1.6.

It follows from (59) that the function

$$f(t) = \sum_k c_k \cos(\omega_k t), \quad c_k \geq 0, \quad (65)$$

is positive definite for any frequencies $\omega_k$ as long as the coefficients $c_k$ are nonnegative.

The coefficients $c_k$ with which the different frequencies contribute to the function $f(t)$ will be very important in what follows and the generic notation $c_k$ will not be adequate for them. We need some notation that incorporates the name of the function $f$, the frequency $\omega_k$, and the fact that we expand $f$ in cosines and not some other periodic functions, specifically not in sines. A popular choice, satisfying all of this criteria is to replace $c_k$ by $\widehat{f^c}(\omega_k)$. So, we write

$$f(t) = \sum \widehat{f^c}(\omega_k) \cos(\omega_k t), \quad \widehat{f^c}(\omega_k) \geq 0, \quad (66)$$

A classical theorem of Bochner says that, conversely, every positive definite function is a limit of functions of the form (66). See Appendix B for more on this.

### 3.1.7.

In equation (64), it is perfectly OK to make the argument of $f$ be a vector $\boldsymbol{t} \in \mathbb{R}^d$. We say that a function $f(\boldsymbol{t})$ is even if

$$f(\boldsymbol{t}) = f(-\boldsymbol{t}),$$

and we say it is *positive definite* if

$$B_f(\boldsymbol{x}) = \sum_{i,j=1}^{n} f(\boldsymbol{t}_i - \boldsymbol{t}_j) x_i x_j \geq 0 \qquad (67)$$

for any $\boldsymbol{x}$ and any $\boldsymbol{t}_1, \ldots, \boldsymbol{t}_n$.

The only modification required in the formula (66) is that the frequencies also become vectors $\boldsymbol{\omega}_k$ and we replace the product $\omega_k t$ by the inner product $(\boldsymbol{\omega}_k, \boldsymbol{t})$. In sum, the function

$$f(t) = \sum \widehat{f^c}(\boldsymbol{\omega}_k) \cos((\boldsymbol{\omega}_k, \boldsymbol{t})), \quad \widehat{f^c}(\boldsymbol{\omega}_k) \geq 0, \qquad (68)$$

is positive definite and every positive definite function of $\boldsymbol{t} \in \mathbb{R}^d$ is a limit of functions of the form (68).

### 3.1.8.

To feed the reader's curiosity, we note briefly that the differences $\boldsymbol{t}_i - \boldsymbol{t}_j$ in the definition of a positive-definite function may be replaced by the ratios $\boldsymbol{t}_i \boldsymbol{t}_j^{-1}$ of elements $\boldsymbol{t}_i$ of an arbitrary group $G$. For the additive group of $\mathbb{R}^d$ we get the positive definite functions as discussed above.

Bochner's theorem is then interpreted as saying that $f$ is a diagonal matrix element of an orthogonal representation of $G$. Viewed from the correct angle, this is very close to a tautology, as noted by Gelfand and Naimark [26] and Segal [48]. See Appendix B for more on this.

### 3.2. The fundamental bound
### 3.2.1.

We will now explain how positive definite functions may be used to bound the density of sphere packings following H. Cohn and N. Elkies [10]. Related considerations, in which translations are replaced by rotations, were used to bound the contact numbers; see [1, 38–40, 42]. As we already mentioned bounds of this type originated in coding theory [16].

### 3.2.2.

Consider an packing of spheres of radius $r$ in $\mathbb{R}^d$ and suppose it is *periodic* in each coordinate direction with period $L$. For instance, consider Figure (6) and let the square in (6) be the square $[0, L]^2 \in \mathbb{R}^2$. Then the leftmost packing in is periodic, and the other two can be made periodic if we erase the spheres intersecting the boundary of the square $[0, L]^2$. In (69) one can see the what the result looks like for the middle packing in (6)

For large $L$, the number of spheres intersecting the boundary of $[0, L]^2$ can be bounded from above by a constant multiple of $L$. Therefore, erasing these spheres changes the density by at most a constant multiple of $L^{-1}$. We conclude that we can come arbitrarily close to the optimal packing density using periodic packings.

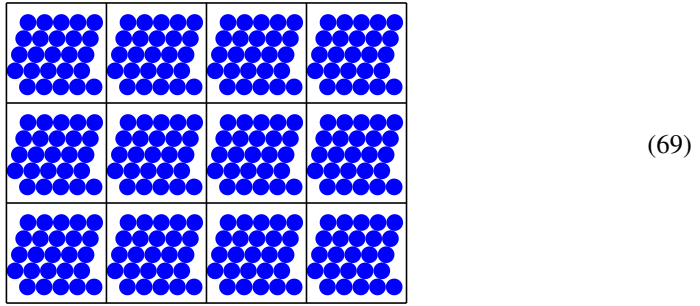

(69)

In dimension $d$, we may need to erase at most a constant multiple of $L^{d-1}$ many spheres, and again this changes the density by at most a constant multiple of $L^{-1}$. We conclude that any upper bound on the density of periodic packings with an arbitrarily large period $L$ gives an upper bound on the density of all sphere packings.

### 3.2.3.

So, returning to our periodic packing, suppose it has $n$ spheres with centers in $[0, L]^d$ and let

$$t_1, \ldots, t_n \in [0, L]^d \tag{70}$$

be the centers of these spheres.

The translates of the cube $[0, L]^d$ tile the whole space $\mathbb{R}^d$ and any periodic sphere packing is just repeated in all these translates. Mathematicians call the basic tile $[0, L]^d$ the *fundamental domain*. So, the number $n$ is the number of spheres per fundamental domain. It directly measures the density of packing by

$$\text{packing density} = \frac{n}{L^d} \operatorname{volume}(B(0, r)), \tag{71}$$

because the sphere packing is periodic. So, our goal is to bound the ratio $n/L^d$.

### 3.2.4.

To bound $n/L^d$, we will use a certain positive definite function $f(t)$, which will be similarly periodic with period $L$ in all coordinates.

To make (68) periodic, the frequencies $\omega$'s should be integer multiples of $\frac{2\pi}{L}$. We define

$$\omega_{\boldsymbol{k}} = \frac{2\pi}{L}\boldsymbol{k}, \tag{72}$$

where $\boldsymbol{k} = (k_1, \ldots, k_d) \in \mathbb{Z}^d$ is a vector with integer entries, and consider a function of the form

$$f(\boldsymbol{t}) = \sum_{\boldsymbol{k}=(k_1,\ldots,k_d)\in\mathbb{Z}^d} \widehat{f^c}(\boldsymbol{k}) \cos\left(\frac{2\pi}{L}(\boldsymbol{k}, \boldsymbol{t})\right), \quad \widehat{f^c}(\boldsymbol{k}) \ge 0. \tag{73}$$

We have written (73) as an infinite sum over all possible frequencies that produce functions with period $L$. Readers who are not comfortable with infinite sums yet may assume that only finitely many of the coefficients $\widehat{f^c}(\boldsymbol{k})$ are nonvanishing in (73). Readers who have seen infinite series, should assume that (73) converges for those values of the argument that will be used below.

The series (73) is a *Fourier series*; see Appendix C for more on this.

### 3.2.5.

Recall that $f(0) > 0$ for any nonzero positive-definite function. For other values of the argument, $f(t)$ may be positive or negative, as exemplified by $\cos \omega t$.

By periodicity, $f$ is positive at any point whose coordinates are integer multiples of $L$. We will denote the set of all such points by $L\mathbb{Z}^d$. Imagine that we managed to arrange $f$ so that

$$\text{distance}(\boldsymbol{t}, L\mathbb{Z}^d) \geq 2r \quad \Rightarrow \quad f(\boldsymbol{t}) \leq 0. \tag{74}$$

In other words, we would like the function $f(t)$ to look like function in Figure (75), namely, positive near the point in $L\mathbb{Z}^d$ and negative away from them:

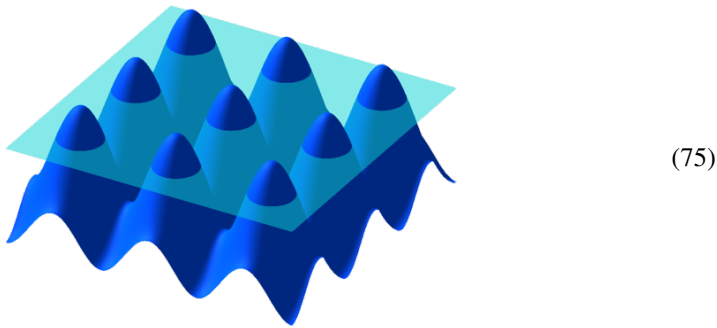

(75)

At this point, this is just a wish. There is absolutely no guarantee that we can find a suitable function. We are only saying that *any* positive definite function satisfying (74) will give some upper bound on the packing density. Whether this bound will be good or bad remains to be seen.

Stressing this logical point is important, because the incredible brilliance of Viazovska's paper [54] is precisely in finding a certain magic positive-definite function that makes everything work.

### 3.2.6.

Because the points (70) are the centers of a periodic sphere packing, we have

$$\text{distance}(\boldsymbol{t}_i - \boldsymbol{t}_j, L\mathbb{Z}^d) \geq 2r, \quad i \neq j,$$

and hence $f(\boldsymbol{t}_i - \boldsymbol{t}_j) \leq 0$ for $i \neq j$. Therefore for the value of (67) at the point $\boldsymbol{x} = (1, \ldots, 1)$ we obtain

$$B_f((1,1,\ldots,1)) = \sum_{i,j} f(\boldsymbol{t}_i - \boldsymbol{t}_j) \leq \sum_i f(\boldsymbol{t}_i - \boldsymbol{t}_i) = nf(0). \tag{76}$$

This will be one side of the eventual inequality involving the number $n$ of spheres in the packing.

### 3.2.7.

For the other side of the inequality, we note from (73) that the matrices of these quadratic forms satisfy

$$B_f = \sum_k \widehat{f^c}(\boldsymbol{k}) \, B_{\cos\left(\frac{2\pi}{L}(\boldsymbol{k},\boldsymbol{t})\right)} \tag{77}$$

and that all terms in this sum are positive definite. Therefore the sum is at least as large as the $\boldsymbol{k} = 0$, that is $B_1$ term

$$B_f((1, 1, \ldots, 1)) \geq \widehat{f}(0) \, B_1((1, 1, \ldots, 1)) = \widehat{f}(0) \, n^2 \,. \tag{78}$$

Here we dropped the superscript from $\widehat{f^c}(0)$ because zero frequency means a constant function and there is no choice between the cosine and sine for $\boldsymbol{k} = 0$. Comparing (78) with (76), we deduce

$$n \leq \frac{f(0)}{\widehat{f}(0)} \,. \tag{79}$$

This is the sought upper bound on the number $n$.

### 3.2.8.

To denominator $\widehat{f}(0)$ in (79) may be interpreted as the *average* of the function $f$ over the fundamental domain $[0, L]^d$. Indeed, the function $\cos\left(\frac{2\pi}{L}(\boldsymbol{k}, \boldsymbol{t})\right)$ with $\boldsymbol{k} \neq 0$ changes sign when shifted by $\frac{1}{2}$ of its minimal period due to

$$\cos(x + \pi) = -\cos(x) \,.$$

Therefore, its average over the whole period vanishes. Thus

$$\begin{aligned}\widehat{f}(0) &= \text{average of } f \text{ over } [0, L]^d \\ &= \frac{1}{L^d} \int_{[0,L]^d} f(\boldsymbol{t}) \, d\boldsymbol{t} \,,\end{aligned} \tag{80}$$

where the second line is for readers familiar with integrals. Readers unfamiliar with integrals may want to consider (80) as the definition of the integral in terms of the average values of the function.

### 3.2.9.

Putting (79) and (80) together, we get

$$\frac{n}{L^d} = \frac{\text{number of spheres in } [0, L]^d}{\text{volume of } [0, L]^d} \leq \frac{f(0)}{\int_{[0,L]^d} f(\boldsymbol{t}) \, d\boldsymbol{t}} \,. \tag{81}$$

Here $f(\boldsymbol{t})$ is a periodic function with period $L$ in each coordinate, which is positive definite and satisfies (74).

### 3.2.10.

As we discussed before, to go from periodic packings to all packings the period $L$ in (81) should get arbitrarily large. Remarkably, there is a way to make *one* function $f$ work for all periods $L$ as follows.

We consider a function $f(t)$ such that

(i) $f(t)$ is positive definite,

(ii) $f(t) \leq 0$ if $\|t\| \geq 2r$,

(iii) $|f(t)|$ decays sufficiently fast as $t \to \infty$.

As before, let (70) be the centers of the spheres in $[0, L]^d$. This means that all centers of the spheres have the coordinates $\{t_i + L\mathbb{Z}^d\}$. If $|f(t)|$ decays sufficiently fast as $t \to \infty$ then the following series

$$\overline{f}(t) = \sum_{v \in L\mathbb{Z}^d} f(t+v) \tag{82}$$

converges, is periodic in $t$, and also positive definite. Evidently,

$$\int_{[0,L]^d} \overline{f}(t)\, dt = \int_{\mathbb{R}^d} f(t)\, dt. \tag{83}$$

As before, we have

$$f(0) \geq \frac{1}{n} \sum_{v \in L\mathbb{Z}^d} \sum_{i,j} f(v + t_i - t_j) = \frac{1}{n} \sum_{i,j} \overline{f}(v + t_i - t_j) \geq \frac{n}{L^d} \int_{\mathbb{R}^d} f(t)\, dt. \tag{84}$$

The first inequality here relies on the fact that $v + t_i - t_j$ is a difference between two sphere centers, and hence has the norm at least $2r$ when nonzero. The second inequality is the inequality (78) applied to the periodic function $\overline{f}$.

We conclude

$$\boxed{\text{density of sphere centers} \leq \min_f \frac{f(0)}{\int_{\mathbb{R}^d} f(t)\, dt}}, \tag{85}$$

where the minimum is over all nonzero functions satisfying the properties (i)–(iii) above.

Bounds of this type are often called *linear programming* bounds because they ask for a minimum of a ratio of two linear functions on a convex set defined by conditions (i)–(iii). Instead of minimizing the ratio, we can consider positive-definite functions normalized by $f(0) = 1$, which is an affine linear equation, and maximize the linear function $\int_{\mathbb{R}^d} f(t)\, dt$ on the resulting convex set.

### 3.2.11.

Note that both the sets and the functions to be extremized are invariant under rotations of $\mathbb{R}^d$, which is a compact group; see Section B.8. Compactness implies there is a well-defined average over all rotations of $f$, which is also a minimizer in (85). This average is a rotation-invariant function, that is, it depends on $\|t\|$ only. Such functions are often called *radial*. To summarize, all we need is a function of one (radial) variable, not $d$ variables.

## 4. Viazovska's magic function
### 4.1. Lattice packings that saturate the bound
#### 4.1.1.

Suppose there are a lattice $\Lambda$ and a function $f$ such that the corresponding packing saturates the bound (85). This implies at once that this packing is the densest possible, but also implies certain very special properties of the function $f$.

Indeed, the inequality (84) was obtained by discarding some nonpositive and nonnegative terms, respectively. If the resulting inequality is an equality then this means all discarded terms vanish.

The first inequality in (84) is an equality if and only if

$$f(v) = 0, \quad \text{for all } v \in \Lambda \setminus \{0\}. \tag{86}$$

If $f$ is radial then it vanishes for all vectors that have the same length as a nonzero vector from $\Lambda$. For $E_8$ this is the set $\sqrt{2n}$, for $n = 1, 2, \ldots$.

#### 4.1.2.

There is a very nice space of functions on $\mathbb{R}^d$ formed by functions that rapidly decay at infinity together with all derivaties. It is called the Schwartz space. For functions $f$ in the Schwarz space, the Fourier transform formulas (161) and (162) from Appendix C become nicely convergent intergrals. The function $\widehat{f}(k) \geq 0$ in

$$f(t) = \int_{\mathbb{R}^d} \widehat{f}(k)\, e^{2\pi i(k,t)}\, dk, \tag{87}$$

is also in Schwarz space and is even/radial if and only if $f(t)$ is even/radial. It is nonnegative because the function $f(t)$ is positive definite by our assumption.

As we will see momentarily, the second inequality in (84) becomes an equality precisely when the Fourier transform vanishes

$$\widehat{f}(k) = 0 \quad \text{for all } k \in \Lambda^\vee \setminus \{0\} \tag{88}$$

for all nonzero vectors in the *dual* lattice, see Appendix C.7. Note the symmetry between (86) and (88). The symmetry is particularly pronounced for $E_8$ because $E_8^\vee = E_8$. For $\Lambda = E_8$ and a radial function $f$, this means the vanishing of the Fourier transform $\widehat{f}(k)$ for all vectors $k$ of length $\sqrt{2n}$, where $n = 1, 2, \ldots$.

#### 4.1.3.

To see (88), let us replace $[0, L]^d$ in the derivation of (84) by the fundamental parallelepiped for $\Lambda$. We redefine

$$\overline{f}(t) = \sum_{v \in \Lambda} f(t+v). \tag{89}$$

Since it is $\Lambda$-periodic, we have

$$\overline{f}(t) = \frac{1}{\sqrt{\Delta_\Lambda}} \sum_{k \in \Lambda^\vee} \widehat{f}(k) \exp(2\pi i(k,t)), \quad \widehat{f}(k) \geq 0, \tag{90}$$

where the coefficients are found from (160) and (162) using

$$\int_{\mathbb{R}^d/\Lambda} \overline{f}(t)\, e^{-2\pi i(k,t)}\, dt = \int_{\mathbb{R}^d} f(t)\, e^{-2\pi i(k,t)}\, dt. \qquad (91)$$

Because there is only *one* sphere in the fundamental parallelepiped, the inequality in (84) becomes

$$\overline{f}(0) = \frac{1}{\sqrt{\Delta_\Lambda}} \sum_{k \in \Lambda^\vee} \widehat{f}(k) \geq \frac{1}{\sqrt{\Delta_\Lambda}} \widehat{f}(0), \qquad (92)$$

where $\frac{1}{\sqrt{\Delta_\Lambda}}$ is the density of the sphere centers, see (46). Clearly, (92) is an equality if and only if (88) holds.

### 4.1.4.

We conclude that to "finish" the proof of optimality of the $E_8$ lattice, one needs to find a function $f(x)$ one variable satisfying the following constraints. We interpret $f(x)$ as a radial function on $\mathbb{R}^8$ and define the Fourier transformed radial function $\widehat{f}(x)$ by

$$\widehat{f}(x) = \int_{\mathbb{R}^8} f(\|t\|) e^{-2\pi i t_1 x}\, dt, \qquad (93)$$

where $t_1$ is the first coordinate of the vector $t$. We need the function $f$ and $\widehat{f}$ to look like the functions in Figure (94)

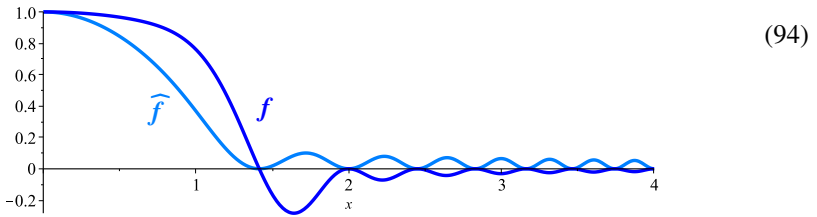

(94)

Namely, function $\widehat{f}(x) \geq 0$ is nonnegative for all $x$ while $f(x) \leq 0$ for $x \geq \sqrt{2}$. Further, both functions vanish for $x = \sqrt{2n}$, $n = 1, 2, \ldots$. Finally, since $\Delta_{E_8} = 1$, we may assume that $f(0) = \widehat{f}(0) = 1$.

### 4.2. The wait is over

In his Fields medal laudatio [8] for Maryna Viazovska, Henry Cohn talks about his attempts to complete this last step, that is, to find the magic function $f$. In particular, he says:

"When Elkies and I proposed this method in 1999, Viazovska was still in secondary school. Without realizing how profoundly difficult the remaining step was, I imagined that we had almost solved the sphere packing problem in eight and twenty-four dimensions, and our inability to find the magic functions was extremely frustrating. At first, I worried that someone else would find an easy solution and leave me feeling foolish for not doing it myself. Over time I became convinced that obtaining these functions was in fact difficult, and others also reached the same conclusion. For example, Thomas Hales has said that *I felt that it would take a Ramanujan to find it* [31]. Eventually, instead of worrying that someone else would solve it, I began to fear that nobody would solve it, and that I would someday die without knowing the outcome. I am grateful that Viazovska found such a satisfying and beautiful solution, and that she introduced wonderful new ideas for the mathematical community to explore."

Viazovska's solution is truly striking. She gives an extremely nontrivial explicit formula for the magic functions in terms of *modular forms*; see Appendix D. There is no way to tell if Ramanujan could have found it, but I would guess that seeing the solution would have made Ramanujan extremely, extremely happy.

Henry Cohn's laudatio [8] contains a very detailed masterfully written account of Viazovska's construction. I do hope the reader feels sufficiently prepared to work through it. While certainly not easy, it is very rewarding. There is a good reason computations like this are recognized by the highest honor in all of mathematics. I also hope the reader opens the interview [55] in which Maryna Viazovska talks, in particular, about her search for the elusive magic function.

### 4.3. Interpolation
#### 4.3.1.

During her search for the magic function, Viazovska conjectured the following systematic way to construct functions like the ones in Figure (94). Namely, she conjectured that a radial Schwartz function on $\mathbb{R}^8$ is uniquely specified by the values of $f$, $f'$, $\widehat{f}$, and $\widehat{f}'$ at the points

$$x = \sqrt{2n}, \quad n = 1, 2, 3 \ldots .$$

In other words, there exists an interpolation basis $a_n$, $b_n$, $\widehat{a}_n$, $\widehat{b}_n$ of the Schwartz space such that for every $f$ we have

$$\begin{aligned} f(x) = &\sum_{n=1}^{\infty} f(\sqrt{2n}) \, a_n(x) + \sum_{n=1}^{\infty} f'(\sqrt{2n}) \, b_n(x) \\ &+ \sum_{n=1}^{\infty} \widehat{f}(\sqrt{2n}) \, \widehat{a}_n(x) + \sum_{n=1}^{\infty} \widehat{f}'(\sqrt{2n}) \, \widehat{b}_n(x) \,. \end{aligned} \quad (95)$$

In particular, the magic function $f$ has to be proportional to $b_1(x)$ because all other coefficients vanish for it.

#### 4.3.2.

This conjecture of Viazovska was proven in her work [12] with Henry Cohn, Abhinav Kumar, Stephen D. Miller, and Danylo Radchenko. They reformulate (95) as a certain functional equation for the following generating series $F(\tau, x)$ and $\widehat{F}(\tau, x)$. By definition,

$$F(\tau, x) = \sum_{n=1}^{\infty} a_n(x) \, e^{2\pi i n \tau} + 2\pi i \tau \sum_{n=1}^{\infty} \sqrt{2n} \, b_n(x) \, e^{2\pi i n \tau}, \quad (96)$$

and similarly for $\widehat{F}(\tau, x)$ with hats everywhere. Note that

$$F(\tau + 2, x) - 2F(\tau + 1, x) + F(\tau, x) = 0, \quad (97)$$

and similarly for $\widehat{F}(\tau, x)$.

The radial function $f_\tau(t) = e^{\pi i \tau \|t\|^2}$, where $t \in \mathbb{R}^8$, has Fourier transform
$$\widehat{f_\tau}(t) = \tau^{-4} e^{-\pi i \|t\|^2/\tau}.$$
Therefore, for $f = f_\tau$ the equation (95) reads
$$e^{\pi i \tau \|t\|^2} = F(\tau, x) + \tau^{-4} \widehat{F}(-1/\tau, x). \tag{98}$$
The authors of [12] solve the equations (98) and (97) in term of modular forms and deduce formulas for the interpolation basis in (96). In particular, this yields a formula for $b_1$, and hence for the $E_8$ magic function.

The appearance of $\tau$ and $-1/\tau$ in the equation (98) is certainly a hint that modular forms have a role to play; compare with Section D.4. Note, however, that $F(\tau, x)$ is *not* periodic in $\tau$, instead (97) says that $(T-1)^2$ annihilates $F(\tau, x)$, where $T$ shifts $\tau$ by 1. (This is a fancy way to say that $F(\tau + n, x)$ is linear in $n$ for $n \in \mathbb{Z}$.) Ultimately this is linked to the appearance of the modular functions for the subgroup $\Gamma(2)$ and also of the quasimodular Eisenstein series $\mathscr{E}_2$.

Like Viazovska's original construction, proving the interpolation formula (95) requires a certain cooperation between math and humans. Math has to make sure there is a miracle to be discovered. Humans have to send their brightest minds on the voyage to discover it.

### 4.3.3.

Similar results are also obtained in [12] for the Leech lattice. These stronger results imply the optimality of $E_8$ and $\Lambda_{24}$ not just for sphere packing but also for certain more general geometric optimizations problems.

### 4.3.4.

I hope the readers share the narrator's sense of awe at this absolutely amazing mathematics and join me in warmest congratulations on it being recognized by the Fields Medal. I also hope the readers got the sense that today's mathematics is not just extraordinarily powerful, but also concrete, understandable, and fun, once one finds the right idea and the right point of view. While finding that right point of view is not at all easy, my biggest hope is to have inspired my youngest readers to believe that mathematics can be beautiful and rewarding, both as a subject and as a profession. Maybe this is also a good place for me to thank Maryna Viazovska and Henry Cohn for this special opportunity to be introduced to their wonderful subject.

## 5. Further reading

The *Quanta Magazine* has published several popular accounts of these and related developments, see [29, 31, 32].

Among introductory or survey articles written by top experts in the field, one could mention [6, 14, 19, 20, 49]. These were written prior to Viazovska's breakthrough. See [7, 8, 36] for expositions of Viazovska's breakthrough.

The reader will surely enjoy reading the textbooks [18, 52] and the comprehensive reference book [15]. A very interesting physics perspective on sphere packings may be found in [45].

I hope the reader has a lot of fun studying these sources as well as the original articles [10–12, 54].

## A. Inner products
### A.1.

In the following discussion we assume that the reader is familiar with basic linear algebra, in particular with the notion of a vector space such as $\mathbb{R}^d$. There exist many beautiful engaging professional expositions of the subject, see for instance [3, 35, 51]. Some readers may find the brief introduction in [44] usable.

A distance function like (1) is an *extra* structure on the linear space $\mathbb{R}^d$, meaning it is not part of the definition of a linear space. But is interacts very nicely with the linear space structures.

First, it is invariant under the translations. So it enough to specify the distance $\|x\|$ to the point $x$ from the origin $0 \in \mathbb{R}^d$. This is also called the *norm* of the vector $x$. The formula

$$\|x\|^2 = \sum x_i^2 \tag{99}$$

is valid in the coordinates with respect to the standard basis $e_1, \ldots, e_d$ of $\mathbb{R}^d$, but will not remain valid in a different basis $e'_1, \ldots, e'_d$. To describe the effect of a linear change of variables, and for many other computations, it is very convenient to introduce the *inner product* associated to (99). By definition

$$(x, y) = \sum_i x_i y_i. \tag{100}$$

The norm (99) and the inner product (100) determine each other by $\|x\|^2 = (x, x)$ and

$$2(x, y) = \|x + y\|^2 - \|x\|^2 - \|y\|^2. \tag{101}$$

If

$$x = \sum x_i e_i = \sum x'_i e'_i \tag{102}$$

is the expansion of $x$ in two different bases, then

$$\|x\|^2 = \sum x_i^2 = \sum_{ij} (e'_i, e'_j) x'_i x'_j. \tag{103}$$

We see that norm squared is given by a *quadratic form* as in Section 3.1.2. Further, this quadratic form is positive definite because $\|x\|^2 > 0$ for any $x \neq 0$.

### A.2.

Given any positive definite quadratic form $B(x)$ we can define a new norm by

$$\|x\|_B^2 = B(x). \tag{104}$$

Using a version of row reduction for the matrix $(b_{ij})$ called the Gram-Schmidt orthogonalization, we can always find a new basis $e'_i$ in which

$$\|x\|_B^2 = \sum (x'_i)^2.$$

Such basis is called an *orthonormal* basis for the form (104).

Linear transformations $\mathsf{g}$ that preserve $\|x\|^2$ are called *orthogonal*. Linear isometries is another word for orthogonal transformations. They take orthonormal bases to orthonormal bases. Writing this condition in terms of matrix entries of $\mathsf{g} = (\mathsf{g}_{ij})$, we see it is equivalent to $\mathsf{g}^T \mathsf{g} = \mathbb{1}$, where $\mathsf{g}^T = (\mathsf{g}_{ji})$ is the transposed matrix and $\mathbb{1}$ denotes the identity matrix. Equivalently

$$\mathsf{g}^{-1} = \mathsf{g}^T, \tag{105}$$

where $\mathsf{g}^{-1}$ is the inverse matrix.

To summarize, invertible linear transformations $\mathsf{g}$ take the standard norm $\|x\|^2$ to all posible positive definite norms (104), and $\mathsf{g}$ takes $\|x\|^2$ to itself if and only if $\mathsf{g}$ is orthogonal. This means that all possible positive definite quadratic forms are the same as invertible linear transformations considered up to precomposing with an orthogonal transformation. See Appendix D for more on this.

### A.3.

If $e_1, \ldots, e_d$ is a basis such that $(e_i, e_j) = 0$ for $i \ne j$ then the expansion $x = \sum x_i e_i$ can be written as

$$x = \sum_i \frac{(x, e_i)}{(e_i, e_i)} e_i. \tag{106}$$

We will find it convenient in Section C.4 below.

Also note the link with the formula for the reflection in the hyperplane orthogonal to the vector $e_1$:

$$r_{e_1}(x) = x - 2 \frac{(x, e_1)}{(e_1, e_1)} e_1. \tag{107}$$

Indeed, the transformation (107) changes the sign of the $e_1$-coefficient in (106) and leaves all other coefficients unchanged. If $(e_1, e_1) = 2$, in particular, if $e_1$ is a root in an even lattice, then (107) simplifies to (38).

### A.4.

Let $\Lambda \subset \mathbb{R}^d$ be a lattice generated by vectors $v_1, \ldots, v_d$ as in Section 1.4. The matrix

$$B_\Lambda = ((v_i, v_j)) \tag{108}$$

is called the Gram matrix. This is a positive definite symmetric matrix and any two collections of vectors $v_1, \ldots, v_d$ and $v'_1, \ldots, v'_d$ with the same Gram matrices (108) can be taken one to another by an orthogonal transformation of $\mathbb{R}^d$. Therefore, in the context of lattice sphere packings, we only care about the Gram matrices of lattices.

### A.5.

Complex numbers are expressions of the form[10]

$$z = a + bi \qquad (109)$$

where $a$ and $b$ are real numbers and $i$ is a symbol satisfying $i^2 = -1$. One can add and multiply complex numbers using this rule. We denote the set of complex numbers by $\mathbb{C}$.

The complex conjugate number is defined by

$$\bar{z} = a - bi \,. \qquad (110)$$

Importantly,

$$z\bar{z} = a^2 + b^2 \,, \qquad (111)$$

which is a nonzero real number for $z \ne 0$. This means, in particular, that $z^{-1} = \frac{1}{a^2+b^2}\,\bar{z}$, which defines division by a nonzero complex number. In other words, complex numbers form a field.

Tuples $z = (z_1, \ldots, z_d)$ of complex numbers form a linear space $\mathbb{C}^d$, in which one defines

$$\|z\|^2 = \sum z_i \bar{z}_i \,, \quad (z, z') = \sum z_i \bar{z}'_i \,. \qquad (112)$$

These norms and inner products are called *Hermitian* and linear transformations that preserve them are called *unitary*.

## B. Groups and positive definite functions
### B.1.

In $\mathbb{R}^3$, consider rotations $\mathsf{g}$ around the origin. For a vector $v \in \mathbb{R}^3$, we will denote by $\mathsf{g}v \in \mathbb{R}^3$ the result of applying the rotation $\mathsf{g}$ to $v$.

Remarkably, if we perform two rotations $\mathsf{g}_2$ and $\mathsf{g}_1$ in succession, the result is another rotation (can you prove this?) which we will denote by $\mathsf{g}_1\mathsf{g}_2$. It is called the composition or the *product* of two rotations. By construction,

$$(\mathsf{g}_1\mathsf{g}_2)v = \mathsf{g}_1(\mathsf{g}_2 v) \,, \qquad (113)$$

for every $v$. Note the order in which we write the product. It is important. In general, $\mathsf{g}_1\mathsf{g}_2 \ne \mathsf{g}_2\mathsf{g}_1$, as we invite the reader to check in in examples. From (113) it follows that

$$(\mathsf{g}_1\mathsf{g}_2)\mathsf{g}_3 = \mathsf{g}_1(\mathsf{g}_2\mathsf{g}_3) \,, \qquad (114)$$

so we do not need the brackets when we write the products.

There is a special identity rotation $\mathbb{1}$ that does nothing and satisfies

$$\mathbb{1}\mathsf{g} = \mathsf{g}\mathbb{1} = \mathsf{g} \,, \qquad (115)$$

---

[10] The numbers $a$ and $b$ are called the real and the imaginary part of the complex number $z$. Terminology notwithstanding, complex numbers really exist. For example, the imaginary unit $i$ is the very first symbol in the Schrödinger equation, one of the fundamental equations describing our really complex world.

for every g. Finally, for every rotation g there is the inverse rotation $g^{-1}$ such that

$$g^{-1}g = gg^{-1} = \mathbb{1}. \qquad (116)$$

**B.2.**

In mathematics, any set G with a special element $\mathbb{1} \in G$, a binary product operation

$$(g_1, g_2) \xrightarrow{\text{product}} g_1 g_2,$$

and a unary inverse operation

$$g \xrightarrow{\text{inverse}} g^{-1},$$

satisfying (114), (115), and (116) is called a group. A subset $G' \subset G$ closed under product and inverse is called a *subgroup*.

This is a very important notion, some examples of which are

$$GL(n, \mathbb{R}) = \text{the group of all invertible } n \times n \text{ real matrices}, \qquad (117)$$
$$O(n, \mathbb{R}) = \text{the subgroup of } n \times n \text{ orhogonal matrices}, \qquad (118)$$
$$SO(n, \mathbb{R}) = \text{orhogonal matrices with det } g = 1, \qquad (119)$$
$$S(n) = \text{permutations of an } n\text{-element set}, \qquad (120)$$

et cetera. Orthogonal matrices were discussed in Section A.2. By construction, we have

$$SO(n, \mathbb{R}) \subset O(n, \mathbb{R}) \subset GL(n, \mathbb{R}) \qquad (121)$$

and we can also embed

$$S(n) \subset O(n, \mathbb{R}) \qquad (122)$$

by making $S(n)$ permute the basis vectors. The group $SO(3, \mathbb{R})$ is the group of rotations of $\mathbb{R}^3$ around the origin discussed above.

For a much much more basic example of a group, one can take the group of real numbers $\mathbb{R}$ with the operation of addition. The zero $0 \in \mathbb{R}$ is the identity element for this operation. Similarly, $\mathbb{R}^d$ is a group with respect to addition. The group $\mathbb{R}^d$ is simpler that the groups in (117) – (120) in one important aspect. The operation in $\mathbb{R}^d$ is *commutative*, meaning that $g_1 g_2 = g_2 g_1$ for any $g_1$ and $g_2$.

In all examples above, $\mathbb{R}$ can be replaced by an arbitrary field. The field $\mathbb{C}$ of complex numbers and the group $U(n)$ of $n \times n$ unitary matrices are particularly important in mathematics.

**B.3.**

One can also use a ring with unit in place of $\mathbb{R}$ above, for instance the ring $\mathbb{Z}$ of integers. In defining $GL(n, \mathbb{Z})$, one needs to make sure that the inverse $g^{-1}$ of an integral matrix $g \in GL(n, \mathbb{Z})$ is also integral. For a commutative ring like $\mathbb{Z}$, it is enough to require that the determinant $\deg g \in \mathbb{Z}$ is an invertible element, meaning that $\deg g = \pm 1$.

The subgroup $GL(n, \mathbb{Z}) \subset GL(n, \mathbb{R})$ consists of matrices that preserve the standard lattice $\mathbb{Z}^n \subset \mathbb{R}^n$. Similarly, matrices preserving an arbitrary lattice $\Lambda \subset \mathbb{R}^n$ form a subgroup that becomes $GL(n, \mathbb{Z})$ in a suitable basis. This is an infinite group. By constrast, orthogonal matrices preserving a given lattice always form a finite group[11]. This group is the boring $\{\pm 1\}$ for a generic lattice, but can be very interesting for lattices like $E_8$ and $\Lambda_{24}$.

**B.4.**

A map between groups

$$\mathsf{G} \to \mathsf{G}'$$

preserving the group structure is called a group homomorphism. A special kind of homomorphism

$$\rho : \mathsf{G} \to O(n, \mathbb{R}) \qquad (123)$$

is called an orthogonal representation of $\mathsf{G}$ of dimension $n$. It represents every $\mathsf{g} \in \mathsf{G}$ by an orthogonal matrix $\rho(\mathsf{g})$ that satisfy

$$\rho(\mathsf{g}_1 \mathsf{g}_2) = \rho(\mathsf{g}_1)\rho(\mathsf{g}_2) \,.$$

For example, (122) is an orthogonal representation. If the target group $O(n, \mathbb{R})$ is replaced by $GL$ or the unitary group, one talks about linear or unitary representation.

Let an orthogonal representation as in (123) be given and and let $v \in \mathbb{R}^n$ be a vector with $\|v\| = 1$. It defines a function on $\mathsf{G}$ by

$$f_{\rho,v}(\mathsf{g}) = (\rho(\mathsf{g})v, v) \,. \qquad (124)$$

Such functions are called diagonal matrix elements. If $v$ is the first basis vector in some basis of $\mathbb{R}^n$ then $f_{\rho,v}(\mathsf{g})$ is the matrix element $\rho(\mathsf{g})_{1,1}$.

Since $\rho(\mathsf{g}^{-1}) = \rho(\mathsf{g})^{-1} = \rho(\mathsf{g})^\mathsf{T}$, we conclude that (124) is symmetric

$$f_{\rho,v}(\mathsf{g}^{-1}) = f_{\rho,v}(\mathsf{g}) \,. \qquad (125)$$

If $\mathsf{g}_1, \ldots, \mathsf{g}_d \in \mathsf{G}$ are arbitrary group elements and $x = (x_1, \ldots, x_n)$ is arbitrary then

$$\left\| \sum_i x_i \rho(\mathsf{g}_i) v \right\|^2 = \sum_{i,j} f_{\rho,v}(\mathsf{g}_i \mathsf{g}_j^{-1}) x_i x_j \,. \qquad (126)$$

Clearly, the quadratic form in (126) is positive semidefinite.

Functions $f(\mathsf{g})$ that are symmetric $f(\mathsf{g}^{-1}) = f(\mathsf{g})$ and produce positive semidefinite forms $\sum_{i,j} f(\mathsf{g}_i \mathsf{g}_j^{-1}) x_i x_j$ are called positive definite functions on $\mathsf{G}$. If $f \neq 0$ then we can normalize it by $f(\mathbb{1}) = 1$.

For the additive group $\mathbb{R}^d$ this is the definition from Section 3.1.5. The analog of Bochner's theorem for $\mathsf{G}$ says that any positive definite function is a diagonal matrix element of an orthogonal representation. This representation could be infinite-dimensional, hence the need for limits in Bochner's theorem. A solid amount of mathematical care is required to work with infinite-dimensional representations, much beyond the introductory style of these notes. We will therefore consider the case of a finite group $\mathsf{G}$, which already contains many key features of the general story.

---

11  Can you prove it? Note that we consider transformations that preserve the origin of $\Lambda$.

**B.5.**

The simplest finite group is the group $\mathbb{Z}/m\mathbb{Z}$ of integers modulo $m$, with the addition operation. It is generated by one element 1, not to be confused with the identity $\mathbb{1}$. In this group, $\mathbb{1}$ is the zero element. For brevity, we denote it by $\mathbb{Z}/m$ in what follows.

The study of representations of groups is a generalization of the theory of eigenvalues and eigenvectors of a matrix. From the eigenvectors of the generator $\rho(1)$, one can conclude that any orthogonal representation of $\mathbb{Z}/m$, in a suitable basis, is the sum of $2 \times 2$ matrix blocks

$$\rho(j) = \begin{pmatrix} \cos \frac{2\pi k j}{m} & -\sin \frac{2\pi k j}{m} \\ \sin \frac{2\pi k j}{m} & \cos \frac{2\pi k j}{m} \end{pmatrix}, \quad k = 1, \ldots, m-1,$$

the trivial representation

$$\rho(j) = 1,$$

and the sign representation $\rho(j) = (-1)^j$, which exists for even $m$. In all cases, the diagonal matrix elements for $\mathbb{Z}/m$ are nonnegative combinations of the functions

$$f_k(j) = \cos \frac{2\pi k j}{m}, \quad k = 0, \ldots, m-1.$$

Our next goal is to show that these exhaust all positive-definite functions $f$ on $\mathbb{Z}/m$.

**B.6.**

To this end, we consider the representation $\rho_{\text{reg}}$ of $\mathbb{Z}/m$ on an $m$-dimensional vector space with basis $\boldsymbol{\delta}_0, \ldots, \boldsymbol{\delta}_{m-1}$ given by

$$\rho_{\text{reg}}(j) \boldsymbol{\delta}_i = \boldsymbol{\delta}_{i+j \bmod m}.$$

This is called the regular representation. We introduce an inner product on it by

$$(\boldsymbol{\delta}_a, \boldsymbol{\delta}_b) = f(b - a).$$

This is symmetric because $f$ is symmetric, positive semidefinite because $f$ is positive definite, and preserved by the action of $\mathbb{Z}/m$. The vectors of zero norm form a linear subspace that is preserved by $\mathbb{Z}/m$, and the representation of $\mathbb{Z}/m$ on the quotient by this subspace is orthogonal. Finally,

$$f(j) = (\rho_{\text{reg}}(j) \boldsymbol{\delta}_0, \boldsymbol{\delta}_0),$$

and this finishes the proof.

**B.7.**

The above discussion may be further simplified if one uses complex numbers and unitary representations. For a unitary representation $\rho(\mathsf{g})$ we have

$$f_{\rho,\boldsymbol{v}}(\mathsf{g}^{-1}) = \overline{f_{\rho,\boldsymbol{v}}(\mathsf{g})}, \tag{127}$$

and (126) turns into a positive definite Hermitian form. A complex-valued function on a group is called positive-definite if it satisfies these two properties.

For a commutative group like $\mathbb{Z}/m$, unitary representations are sums of 1-dimensional representations

$$\chi_k(j) = \exp\left(\frac{2\pi i j k}{m}\right), \quad k = 0, \ldots, m-1, \tag{128}$$

and the argument given above proves that positive definite functions on $\mathbb{Z}/m$ are nonnegative linear combinations of the functions (128). One-dimensional representations are called *characters* and often denoted by the letter $\chi$.

In (128), the letter $i$ denotes the imaginary unit, and the exponential of an imaginary number may be defined by

$$\begin{aligned} e^{it} &= 1 + it + \frac{(it)^2}{2} + \frac{(it)^3}{3!} + \ldots \\ &= \left(1 - \frac{t^2}{2} + \frac{t^4}{4!} - \ldots\right) + i\left(t - \frac{t^3}{3!} + \frac{t^5}{5!} - \ldots\right) \\ &= \cos(t) + i\sin(t), \end{aligned} \tag{129}$$

as discovered by Leonhard Euler around 1740. Note the famous special case

$$e^{\pi i} = -1. \tag{130}$$

### B.8.

We started this section with a discussion of rotations of $\mathbb{R}^3$, which in addition to forming a group have two further important properties. First, rotations form a manifold, namely the real projective 3-space. Groups forming a manifold are called Lie groups. Second, this manifold is compact.

Compact Lie groups are very important in mathematics and they have been completely classified. This classification includes the classification of compact connected Lie groups and of finite simple groups. In both cases, there are certain well-understood infinite series as well as finitely many exceptional cases, surrounded by a much denser air of mystery. The classification of compact connected Lie groups is very close to the ADE classification[12] from Section 2.3.3, and it ends in the largest exceptional group $E_8$. Among the finite groups, there is the largest sporadic group called the Monster, which is very closely connected to the Leech lattice $\Lambda_{24}$.

## C. Fourier series
### C.1.

Let us revisit the regular representation of $\mathbb{Z}/m$ from Section B.6. Every group $\mathsf{G}$ acts on the linear space of functions $f : \mathsf{G} \to \mathbb{C}$ by the following rule:

$$\left[\rho_{\mathrm{reg}}(\mathsf{g})f\right](\mathsf{g}') = f(\mathsf{g}'\mathsf{g}). \tag{131}$$

---

[12] In this classification, the lattices $A_n$ correspond to special unitary groups $SU(n+1)$, while the lattices $D_n$ correspond to even orthogonal groups $SO(2n)$.

In (131) we have the result of evaluation of the new function $\rho_{\text{reg}}(\mathsf{g})f$ at a group element $\mathsf{g}'$. The square brackets in (131) are put around $\rho_{\text{reg}}(\mathsf{g})f$ just to stress that this this a new function, obtained by the action of $\mathsf{g}$ from the original function $f$.

The basis $\delta_0, \ldots, \delta_{m-1}$ from Section B.6 corresponds the functions

$$\delta_k(j) = \delta_{kj}, \text{ where } \delta_{kj} = \begin{cases} 1, & k = j, \\ 0, & \text{otherwise}. \end{cases} \quad (132)$$

### C.2.

The following Hermitian product

$$(f_1, f_2)_{\text{reg}} = \sum_{j \in \mathbb{Z}/m} f_1(j) \overline{f_2(j)} \quad (133)$$

makes the regular representation of $\mathbb{Z}/m$ unitary. The basis (132) satisfies

$$(\delta_k, \delta_{k'})_{\text{reg}} = \delta_{kk'}. \quad (134)$$

### C.3.

Now consider the functions (128) as elements of the regular representation. We compute

$$(\chi_k, \chi_{k'})_{\text{reg}} = \sum_{j \in \mathbb{Z}/m} \exp\left(\frac{2\pi i j (k - k')}{m}\right) = |\mathsf{G}|\, \delta_{k,k'}. \quad (135)$$

Indeed, the sum in (135) is a sum of a geometric progression and it vanishes if $k \neq k'$. We put the cardinality $|\mathsf{G}|$ in (135) instead of $m$ because (135) is the simplest case of a very general relations known as orthogonality of characters (and other matrix elements of irreducible representation). In the case of a finite group, the cardinality $|\mathsf{G}|$ is the correct factor to put into these orthogonality relations.

### C.4.

We have found two orthogonal bases $\{\delta_k\}$ and $\{\chi_k\}$ in the space of complex-valued functions on $\mathsf{G} = \mathbb{Z}/m$. Let us expand a general function in these bases using (106).

The expansion in the basis $\{\delta_k\}$ amounts to a tautology:

$$f = \sum_{k \in \mathbb{Z}/m} f(k)\, \delta_k. \quad (136)$$

The expansion in the basis $\chi_k$, by contrast, amounts to something very nontrivial. By (106), the coefficients $\widehat{f}(k)$ in the expansion

$$f = \sum_{k \in \mathbb{Z}/m} \widehat{f}(k)\, \chi_k \quad (137)$$

are given by

$$\widehat{f}(k) = \frac{(f, \chi_k)_{\text{reg}}}{(\chi_k, \chi_k)_{\text{reg}}} = \frac{1}{m} \sum_{j=0}^{m-1} f(j)\, \exp\left(-\frac{2\pi i j k}{m}\right). \quad (138)$$

The expansion (137), written out, takes a very similar form

$$f(j) = \sum_{k=0}^{m-1} \widehat{f}(k) \exp\left(\frac{2\pi i j k}{m}\right). \tag{139}$$

Formulas (138) and (139) describe the *Fourier transform* on the commutative group $\mathsf{G} = \mathbb{Z}/m$. By (139) every function $f(j)$ can be written as a combination of characters $\chi_k$. The coefficients $\widehat{f}(k)$ in (139) are the average values of $f\overline{\chi_k}$.

A similar Fourier transform on groups exists very generally. For noncommutative groups, one should take matrix elements of unitary representations instead of characters. This will remain entirely outside of our narrative.

### C.5.

After talking about Fourier transform for a finite commutative group $\mathbb{Z}/m$, let us consider the simplest commutative connected Lie group $SO(2)$ of rotations in $\mathbb{R}^2$ around the origin. A rotation is specified by the angle $\phi \in [0, 2\pi]$, with the endpoints $0$ and $2\pi$ representing the same identity element in $SO(2)$. The group operation is the addition of angles $\phi$, taken modulo $2\pi$. Thus, we may think of $SO(2)$ as the quotient

$$SO(2) = \mathbb{R}/2\pi\mathbb{Z},$$

of a linear group by a lattice subgroup.

For any $m$, we have the subgroup

$$\mathbb{Z}/m = \left\{\phi = \frac{2\pi j}{m}, j = 0, \ldots, m-1\right\} \subset SO(2) \tag{140}$$

formed by rotations that preserve a regular $m$-gon. As $m$ gets large, these become denser and denser. Let us *formally* take the $m \to \infty$ limit in the formulas (138) and (139), and see if we get the formulas for the Fourier transform on $SO(2)$.

Let us rewrite the formulas (138) and (139) using the variable $\phi = \frac{2\pi j}{m}$. We get

$$\widehat{f}(k) = \frac{1}{m} \sum_{\phi \in \mathbb{Z}/m} f(\phi) e^{-ik\phi}, \tag{141}$$

$$f(\phi) = \sum_{k \in \mathbb{Z}/m} \widehat{f}(k) e^{ik\phi}. \tag{142}$$

In (141), we interpret $\mathbb{Z}/m$ as the subgroup (140), while in (142) we have a summation of a periodic function of $k$ over any period of length $m$ in $\mathbb{Z}$. As $m \to \infty$, the sum in (141) approximates the integral over the group $SO(2)$, while the sum (142) becomes the sum over all integers $k$. Thus, we get

$$\widehat{f}(k) = \frac{1}{2\pi} \int_0^{2\pi} f(\phi) e^{-ik\phi} d\phi, \tag{143}$$

$$f(\phi) = \sum_{k \in \mathbb{Z}} \widehat{f}(k) e^{ik\phi}. \tag{144}$$

We stress that our derivation of these formulas was just by a formal analogy with the case of a finite group and much more serious work is required to both interpret these formulas correctly and prove them. Questions like these belong to the field of harmonic analysis, which is a very deep and important part of mathematics. We leave the reader by the entrance to this glorious edifice, referring to [46, 50] for possible further reading. But to stimulate the reader's curiosity, we will do one example.

### C.6.

Fix some angle $\phi_0$ and consider the function

$$f(\phi) = \begin{cases} 1, & \cos(\phi) \geq \cos(\phi_0), \\ 0, & \text{otherwise}. \end{cases}$$

In other words, this function equal 1 on the "spherical cap" $[-\phi_0, \phi_0]$ and vanishes outside of it. From (143), we compute

$$\widehat{f}(k) = \frac{1}{2\pi} \int_{-\phi_0}^{\phi_0} e^{-ik\phi} d\phi = \begin{cases} \frac{\phi_0}{\pi}, & k = 0, \\ \frac{\sin(k\phi_0)}{k\pi}, & k \neq 0, \end{cases} \qquad (145)$$

where we used the relation (129) between the complex exponential and the trigonometric functions. Using (129) again, we can write (144) as follows

$$f(\phi) \stackrel{?}{=} \frac{\phi_0}{\pi} + 2 \sum_{k=1}^{\infty} \frac{\sin(k\phi_0)\cos(k\phi)}{k\pi}, \qquad (146)$$

where the question mark indicates that the exact interpretation of this equality is beyond the scope of these notes.

A picture being worth a thousand words, we just plot the partial sums of the series above for $\phi_0 = \frac{\pi}{3}$ and $k$ up to 3, 5, 10, and 50, respectively.

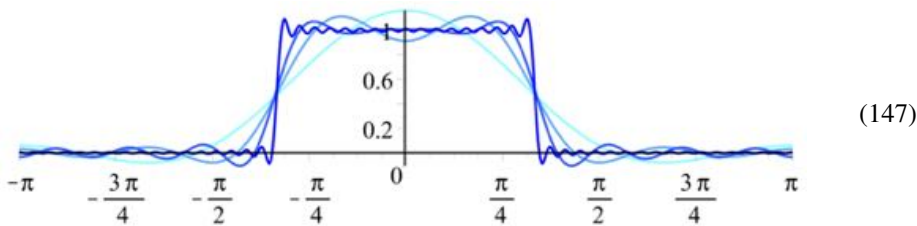

(147)

One salient feature of (147) are the very strong oscillations of the Fourier series near the point of discontinuity of the function, known as the Gibbs phenomenon.

Experimenting with Fourier series is a lot of fun, and we invite the reader to do more experiments! The functions that are actually needed in Viazovska's proof are infinitely differentiable and their Fourier expansions converge to them very nicely.

### C.7.

The following common generalization of (138), (139), (143), (144) is valid for any commutative Lie group[13] G. It describes the expansion of a function $f$ on G in terms of the characters of G.

---

[13] as well as for more general locally compact commutative groups

Unitary characters of $\mathsf{G}$, that is, continuous homomorphisms

$$\chi : \mathsf{G} \to U(1) = \{z \in \mathbb{C}, |z| = 1\}, \tag{148}$$

form a commutative group $\mathsf{G}^\wedge$ with respect to pointwise multiplication of characters. The trivial character $\chi = 1$ is the identity of this group. If the group $\mathsf{G}$ is compact then $\mathsf{G}^\wedge$ is discrete, and visa versa. The group $\mathsf{G}^\wedge$ is called the Pontryagin dual group, or the dual group for short.

Mathematicians write

$$1 \to \mathsf{G}_1 \to \mathsf{G} \to \mathsf{G}_2 \to 1 \tag{149}$$

to indicate that $\mathsf{G}_1$ is a Lie subgroup of $\mathsf{G}$ with quotient $\mathsf{G}_2$. They call a sequence of the form (149) a short exact sequence. Duality reverses short exact sequences:

$$1 \to \mathsf{G}_2^\wedge \to \mathsf{G}^\wedge \to \mathsf{G}_1^\wedge \to 1 \tag{150}$$

which means that characters of $\mathsf{G}_2$ are the characters of $\mathsf{G}$ that are trivial when restricted to $\mathsf{G}_1$ and vice versa. One replaces the 1's by 0's in short exact sequences when the group operation is written as addition.

For example, any inner product $(\cdot, \cdot)$ on $\mathbb{R}^d$ gives the identification $(\mathbb{R}^d)^\wedge \cong \mathbb{R}^d$ by

$$\chi_{\boldsymbol{k}}(\boldsymbol{t}) = \exp(2\pi i (\boldsymbol{k}, \boldsymbol{t})). \tag{151}$$

If $\Lambda \subset \mathbb{R}^d$ is a lattice then the quotient group in

$$0 \to \Lambda \to \mathbb{R}^d \to \mathbb{R}^d/\Lambda \to 0 \tag{152}$$

is a group abstractly isomorphic to $SO(2)^d$ which can be realized concretely by gluing the opposite sides of the fundamental parallelepiped for $\Lambda$. Mathematicians call such group a *torus*. Using (151), we get the identifications

$$\left(\mathbb{R}^d/\Lambda\right)^\wedge = \Lambda^\vee, \tag{153}$$

$$\Lambda^\wedge = \mathbb{R}^d/\Lambda^\vee, \tag{154}$$

where $\Lambda^\vee$ is the *dual lattice*

$$\Lambda^\vee = \{\boldsymbol{k}, \text{ such that } (\boldsymbol{k}, \boldsymbol{v}) \in \mathbb{Z} \text{ for all } \boldsymbol{v} \in \Lambda\}. \tag{155}$$

### C.8.

While as abstract groups, all lattices and tori of the same dimension are isomorphic, they are all very different in the context of sphere packing and other problems involving distances and inner products. It is, therefore, important to distinguish clearly between a lattice $\Lambda$ and the dual lattice $\Lambda^\vee$.

In general, $\Lambda^\vee$ is very different from $\Lambda$. For example, if we scale $\Lambda$ by a factor then $\Lambda^\vee$ scales by the reciprocal factor. However,

$$(\mathbb{Z}^d)^\vee = \mathbb{Z}^d, \quad E_8^\vee = E_8, \tag{156}$$

and, in general, if a lattice is *integral*, which means that $(\boldsymbol{v}_1, \boldsymbol{v}_2) \in \mathbb{Z}$ for all $\boldsymbol{v}_1, \boldsymbol{v}_2 \in \Lambda$, and *unimodular*, which means that $\Delta_\Lambda = 1$ then[14] $\Lambda^\vee = \Lambda$.

---

14    Indeed, integrality implies that $\Lambda \subset \Lambda^\vee$, while $\Delta_\Lambda$ is the order of the group $\Lambda^\vee/\Lambda$.

### C.9.

The Fourier transform on a general *compact commutative* group $\mathsf{G}$ takes the form

$$f(\mathsf{g}) = \sum_{\mathbf{k} \in \mathsf{G}^\wedge} \widehat{f}(\mathbf{k})\, \chi_{\mathbf{k}}(\mathsf{g})\,, \tag{157}$$

$$\widehat{f}(\mathbf{k}) = \int_{\mathsf{G}} f(\mathsf{g})\, \overline{\chi_{\mathbf{k}}(\mathsf{g})}\, d_{\text{prob}}\mathsf{g}\,, \tag{158}$$

where the integration is with respect to the invariant measure on the group $\mathsf{G}$ of total volume 1. Measures of total volume 1 are often called *probability* measure, hence the subscript in (158).

For tori, the integral in (158) is just the usual integral over the fundamental domain, normalized so that the volume of the fundamental domain equals 1. Recall that this volume equals $\sqrt{\Delta_\Lambda}$ with respect to usual volume form $d\mathbf{t}$. Therefore, for $\mathsf{G} = \mathbb{R}^d/\Lambda$, the Fourier transform takes the form

$$f(\mathbf{t}) = \sum_{\mathbf{k} \in \Lambda^\vee} \widehat{f}(\mathbf{k})\, e^{2\pi i (\mathbf{k}, \mathbf{t})}\,, \tag{159}$$

$$\widehat{f}(\mathbf{k}) = \frac{1}{\sqrt{\Delta_\Lambda}} \int_{\mathbb{R}^d/\Lambda} f(\mathbf{t})\, e^{-2\pi i (\mathbf{k}, \mathbf{t})}\, d\mathbf{t}\,. \tag{160}$$

### C.10.

We remind the reader that we glide over all the deep analytic issues involved in the Fourier transform on continuous groups. For our narrative, this is justified by the fact that the actual functions that come up in Viazovska's proof have very nice analytic properties.

While for noncompact groups $\mathsf{G}$ the Fourier transform presents furter analytic difficulties, one can formally take the limit of a very large lattice $\Lambda$ in (159) and (160) and obtain

$$f(\mathbf{t}) = \int_{\mathbb{R}^d} \widehat{f}(\mathbf{k})\, e^{2\pi i (\mathbf{k}, \mathbf{t})}\, d\mathbf{k}\,, \tag{161}$$

$$\widehat{f}(\mathbf{k}) = \int_{\mathbb{R}^d} f(\mathbf{t})\, e^{-2\pi i (\mathbf{k}, \mathbf{t})}\, d\mathbf{t}\,. \tag{162}$$

As $\Lambda$ becomes very large, the dual lattice $\Lambda^\vee$ becomes very dense and the sum in (159) becomes the integral in (161).

## D. Modular forms

### D.1. The space of lattices

The many special lattices we met in these notes may be interpreted as some very special points in a space that parametrizes all possible lattices $\Lambda \subset \mathbb{R}^d$. How should think about this space?

An arbitrary lattice $\Lambda \subset \mathbb{R}^d$ may be obtained from the standard lattice $\mathbb{Z}^d \subset \mathbb{R}^d$ by a change of basis or, equivalently, as a result of linear transformation

$$\Lambda = \mathsf{g}\mathbb{Z}^d\,, \quad \mathsf{g} \in GL(d, \mathbb{R})\,. \tag{163}$$

Further, $g\mathbb{Z}^d = \mathbb{Z}^d$ if and only if $g \in GL(d,\mathbb{Z})$. Thus

$$\{\text{lattices in } \mathbb{R}^d\} = GL(d,\mathbb{R})/GL(d,\mathbb{Z}), \quad (164)$$

where the quotient sign means that we identify $g_1$ and $g_2$ if $g_1^{-1}g_2 \in GL(d,\mathbb{Z})$.

For sphere packing and many other problems, we don't want to distinguish between isometric lattices, that is, lattices that differ by postcomposing $g$ with an orthogonal transformation. Thus we consider

$$\left\{\begin{array}{c}\text{lattices in } \mathbb{R}^d \\ \text{up to isometry}\end{array}\right\} = O(d,\mathbb{R})\backslash GL(d,\mathbb{R})/GL(d,\mathbb{Z}). \quad (165)$$

Finally, for the sphere packing problem, we can rescale the lattice arbitrarily, while simultaneously rescaling the radius of the spheres. Thus, one may want to consider

$$\left\{\begin{array}{c}\text{lattices in } \mathbb{R}^d \text{ up to} \\ \text{scale and isometry}\end{array}\right\} = (\mathbb{R}_{>0}\, O(d,\mathbb{R}))\backslash GL(d,\mathbb{R})/GL(d,\mathbb{Z}), \quad (166)$$

where $\mathbb{R}_{>0}$ is the subgroup of $GL(d,\mathbb{R})$ consisting of positive multiples of the identity matrix.

### D.2.

Let's see what the space (166) looks like for $d = 2$. Let $\Lambda$ be a lattice and let $v \in \Lambda$ be a vector of minimal length. We will complete $v$ to a basis $\{v, v'\}$ of the lattice $\Lambda$ by choosing a shortest vector $v'$ among those not proportional to $v$. Note, however, that $-v'$ is another vector with the same properties as $v'$.

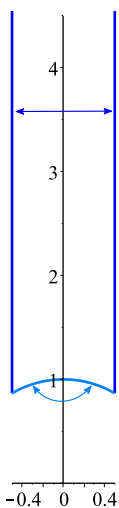

Since we take lattices up to scale and isometry we may arrange so that $v = e_1 = (1,0)$ is the standard basis vector. What are the possibilities for $v' = (v_1', v_2')$? First, we need to have

$$\|v'\|^2 = v_1'^2 + v_2'^2 \geq 1. \quad (167)$$

Second,

$$\|v' \pm e_1\|^2 = \|v'\|^2 \pm 2v_1' + 1, \quad (168)$$

which means that $v'$ can be made shorter by adding or subtracting $e_1$ unless

$$|v_1'| \leq \frac{1}{2}. \quad (169)$$

The combination of (167) and (169) describes the domain shown in the figure on the left plus the symmetric domain below the horizontal axis. Using the symmetry between $v'$ and $-v'$, we may restrict our attention to the figure on the left. Every point in this domain corresponds to a certain lattice in $\mathbb{R}^2$, but some pairs of points on the boundary correspond to the same lattices. Indeed the vertical boundaries differ by a shift by $e_1$. Since this shift does not change the lattice, they have to be glued as indicated. Points on the round boundary correspond to lattices generated by two vectors $v, v'$ of equal length. Since we can declare either one of them to be equal to $e_1$, we may demand additionally that the angle between $v$ and $v'$ is not larger than $\pi/2$. This leads to gluing the round boundary to itself as indicated.

Once we do this gluing, we get a surface of the shape shown in Figure (170):

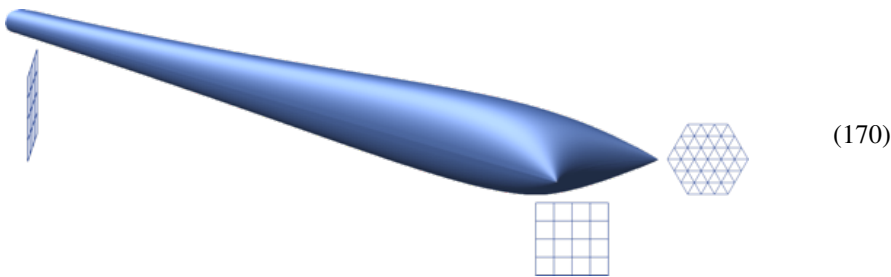

(170)

The two cone points of this surface correspond to the special lattices — the square lattice and the hexagonal latices. They are special because they are preserved by some nontrivial transformations in $O(2,\mathbb{R})$. The surface in (170) has an infinitely long neck, which is sometimes call the *cusp*. A lattice $\Lambda$ runs off to infinity in this neck if $\frac{\|v'\|}{\|v\|} \to \infty$.

### D.3.

Functions on quotient spaces like (164), (165), and (166) are called *automorphic functions*. They are objects of extreme beauty, complexity, and importance for mathematics. Suffices to say that they played a very essential role in Andrew Wiles' proof of Fermat's last theorem.

In these notes, we will limit ourselves to the discussion of one basic class of such functions for $d = 2$, called the Eisenstein series[15]. Let $\Lambda \subset \mathbb{R}^2$ be a lattice. We can identify $\mathbb{R}^2$ with the complex numbers $\mathbb{C}$ from Section A.5 and then $\Lambda$ becomes a subset of $\mathbb{C}$. We define

$$\mathcal{E}_k(\Lambda) = \tfrac{1}{2} \sum_{z \in \Lambda_{\text{primitive}}} \frac{1}{z^k}, \qquad (171)$$

where primitive means that $z$ is not a positive multiple of another vector, in particular this means that $z \neq 0$. The series (171) converges absolutely for $k > 2$ and vanishes for $k$ odd because the contributions of $z$ and $-z$ cancel. For even $k$, $z$ and $-z$ make the same contribution, hence the $\tfrac{1}{2}$ factor in front.

If we multiply the lattice $\Lambda$ by a complex number $w$ then

$$\mathcal{E}_k(w\Lambda) = w^{-k} \mathcal{E}_k(\Lambda). \qquad (172)$$

Note that multiplication by $w$ combines rotations and scaling of $\Lambda$. As a result, Eisenstein series are not exactly invariant under rotation and scaling, but rather transform in the way described by (172) under rotation and scaling. Using (172), we may assume that

$$v = 1, \quad v' = \tau, \qquad (173)$$

where $\tau$ is a complex number in the upper half-plane. The series (171) being a sum of the terms $(n + m\tau)^{-k}$, where $(n, m)$ runs over coprime pairs of integers, the Eisenstein series $\mathcal{E}_k$ is a holomorphic function of the parameter $\tau$.

---

**15**    or, more precisely, holomorphic Eisenstein series

### D.4.

Exchanging the roles of $v$ and $v'$ in (173) leads to the transformation

$$\tau \mapsto -1/\tau, \qquad (174)$$

which takes the complex upper half-plane to itself. Here is what this transformation does to the Cartesian coordinates on the upper half-plane:

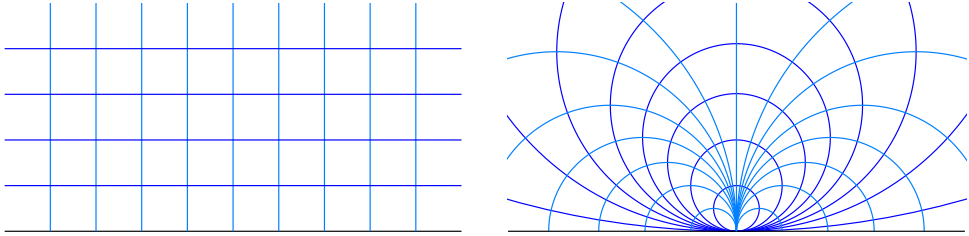

### D.5.

The following beautiful formulas for the series $\mathscr{E}_k$ may be derived in terms of the variable $q = e^{2\pi i \tau}$. For $\tau$ in the upper half-plane, the corresponding $q$ lies in the unit circle $|q| < 1$. We have

$$\mathscr{E}_k = 1 - \frac{2k}{B_k} \sum_{n=1}^{\infty} \left( \sum_{d \mid n} d^{k-1} \right) q^n, \qquad (175)$$

where $B_k$ are the Bernoulli numbers, and the coefficient of $q^n$ is determined by summing over all divisors $d$ of the number $n$. In particular, we have

$$\mathscr{E}_2 = 1 - 24\,q - 72\,q^2 - 96\,q^3 - 168\,q^4 - 144\,q^5 - 288\,q^6 - 192\,q^7 - 360\,q^8 - \ldots, \qquad (176)$$

$$\mathscr{E}_4 = 1 + 240\,q + 2160\,q^2 + 6720\,q^3 + 17520\,q^4 + 30240\,q^5 + 60480\,q^6 + 82560\,q^7 + 140400\,q^8 + \ldots, \qquad (177)$$

$$\mathscr{E}_6 = 1 - 504\,q - 16632\,q^2 - 122976\,q^3 - 532728\,q^4 - 1575504\,q^5 - 4058208\,q^6 - 8471232\,q^7 - 17047800\,q^8 - \ldots, \qquad (178)$$

where we have added the series $\mathscr{E}_2$. Not being absolutely convergent, the series $\mathscr{E}_2$ may be summed with some further choices and (176) is the result. The readers who feel they have already seen the number 240 somewhere recently are not mistaken.

### D.6.

A holomorphic function $f(\Lambda)$ of a lattice $\Lambda \subset \mathbb{C}$ which satisfy (172) and remains bounded as long as the shortest vector $v \in \Lambda$ is bounded away from 0 is called a modular form of weight $k$. We can multiply modular forms of different weights and the weights add under multiplication. Thus modular forms form an algebra, and it is a classical theorem that

$$\text{Modular forms} = \mathbb{C}[\mathscr{E}_4, \mathscr{E}_6]. \qquad (179)$$

The square brackets mean that $\mathscr{E}_4$ and $\mathscr{E}_6$ generate the algebra of modular forms freely, meaning, they do not satisfy any polynomial equation in two variables.

One often adds the series $\mathscr{E}_2$, whose converges requires some regularization making its transformation law a bit more complicated. With this addition, the algebra (179) becomes

$$\text{Quasimodular forms} = \mathbb{C}[\mathscr{E}_2, \mathscr{E}_4, \mathscr{E}_6]. \qquad (180)$$

**D.7.**

There are countless applications of modular forms to the study of lattices. A very major one is the subject of these introductory notes — Viazovska's gigantic breakthrough. For a much more basic one, consider the following situation.

Let $\Lambda \subset \mathbb{R}^d$ be a lattice. We can associate to it its theta series

$$\Theta_\Lambda(q) = \sum_{v \in \Lambda} q^{\frac{1}{2}\|v\|^2}. \tag{181}$$

This converges for $|q| < 1$ for any lattice $\Lambda$. If $\Lambda$ is even then this is series in $q$. And if $\Lambda$ is additionally unimodular than this is modular form of weight $d/2$.

In particular for $\Lambda = E_8$ we should get a modular form of weight 4 and from (179) we see that it can only be a multiple of $\mathscr{E}_4$. Comparing the coefficients of $q^0$, we conclude

$$\Theta_{E_8} = \mathscr{E}_4. \tag{182}$$

Thus the coefficients in (177) count the vectors of a length $\sqrt{2n}$ in the lattice $E_8$. In particular, 240 is the number of roots.

### E. The volume of a $d$-dimensional ball

**E.1.**

Let $B(0, r)$ be the $d$-dimensional ball (2) of radius $r$. Its volume is proportional to $r^d$, namely

$$\text{Vol}\, B(0, r) = \mathsf{v}_d\, r^d, \tag{183}$$

with some proportionality constant $\mathsf{v}_d$. Our goal in this section is to compute this constant. As we will see, it is given in terms of a certain special function (184).

**E.2.**

The Gamma function is defined by the following integral

$$\Gamma(s) = \int_0^\infty e^{-x} x^{s-1} dx, \tag{184}$$

which converges when $s > 0$. For complex $s$, the integral (184) converges when the real part $\Re s > 0$. Integration by parts gives

$$\Gamma(s + 1) = s\Gamma(s), \tag{185}$$

and using this formula one can extend the definition of $\Gamma(s)$ to all values of $s$ except $s = 0, -1, -2, \ldots$.

From (185) and the base case $\Gamma(1) = 1$, we conclude

$$\Gamma(n) = 1 \cdot 2 \cdot 3 \cdots (n-1) = (n-1)!, \quad n = 1, 2, 3, \ldots. \tag{186}$$

The Gamma functions is, in a certain technical sense, the most natural extension of the factorial $(n-1)!$ to a function of a complex variable.

### E.3.

Consider the plot (thick curve) of the logarithm of integrand in (184) for $s = 1000$.

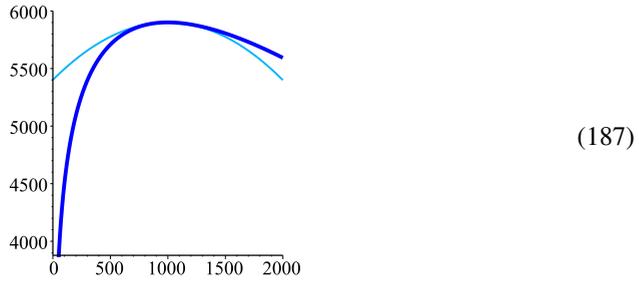

(187)

We will come back to the meaning of the thin curve later. In (187), we have the logarithm, meaning the integrand itself takes very large values. Their maximum is at the solution of

$$(e^{-x}x^{s-1})' = \left(-1 + \frac{s-1}{x}\right)e^{-x}x^{s-1} = 0 \quad \Rightarrow \quad x = s - 1. \tag{188}$$

Hence $\Gamma(1000)$ should be something of the order $\left(\frac{999}{e}\right)^{999}$. Refining this argument, one can deduce a more precise asymptotic relation

$$\Gamma(s+1) \sim \sqrt{2\pi s}\left(\frac{s}{e}\right)^s \tag{189}$$

known as the Stirling formula. It is often used to approximate factorials.

### E.4.

Let's put $x = y^2$ in (184). Since $dx = 2y\,dy$, we get

$$\Gamma(s) = 2\int_0^\infty e^{-y^2} y^{2s-1} dy. \tag{190}$$

In particular, we get the famous Gaussian integral for $s = \frac{1}{2}$

$$\int_{-\infty}^\infty e^{-y^2} dy = \Gamma(\tfrac{1}{2}). \tag{191}$$

Let us multiply $d$ copies of (191). We get

$$\begin{aligned}\Gamma(\tfrac{1}{2})^d &= \int_{\mathbb{R}^d} e^{-(y_1^2 + \cdots + y_d^2)} dy_1 \ldots dy_d \\ &= \int_{\mathbb{R}^d} e^{-\|\mathbf{y}\|^2} d\mathbf{y}\end{aligned} \tag{192}$$

We observe the remarkable fact that the integrand in (192) depends only on the norm of the vector $\mathbf{y}$.

While in this section we have to assume that the reader has some familiarity with integrals, it may be worth recalling how Lebesgue integral are defined. One approximates the integrand by a function taking a discrete set of values and weights each value by the volume of the set where this value is taken.

In particular, we can approximate the function $\|y\|$ by the functions

$$\varepsilon \left\lceil \frac{\|y\|}{\varepsilon} \right\rceil \to \|y\|, \quad \varepsilon \to 0, \tag{193}$$

that take the value $r = k\varepsilon$, $k = 0, 1, 2, \ldots$, on the spherical shell formed by the difference of $B(0, r)$ and the smaller ball $B(0, r - \varepsilon)$. From (183), we conclude

$$\operatorname{Vol} B(0, r) - \operatorname{Vol} B(0, r - \varepsilon) \approx d\, \mathsf{v}_d\, r^{d-1} \varepsilon. \tag{194}$$

Therefore

$$\int_{\mathbb{R}^d} e^{-\|y\|^2}\, dy = d\, \mathsf{v}_d \int_0^\infty e^{-r^2} r^{d-1}\, dr = \tfrac{d}{2} \mathsf{v}_d \Gamma(\tfrac{d}{2}) = \mathsf{v}_d \Gamma(\tfrac{d}{2} + 1), \tag{195}$$

where we have used equalities (190) and (185).

Putting (192) and (195) together, we conclude

$$\mathsf{v}_d = \frac{\Gamma(\tfrac{1}{2})^d}{\Gamma(\tfrac{d}{2} + 1)}. \tag{196}$$

### E.5.

To simplify (196), we note that the $\pi r^2$ formula for the area of circle computes the Gaussian integral! Indeed,

$$\mathsf{v}_2 = \pi \quad \Rightarrow \quad \Gamma(\tfrac{1}{2}) = \sqrt{\pi}. \tag{197}$$

In fact, the $\sqrt{2\pi s}$ prefactor in the Stirling formula (189) comes from approximating the integral (184) by a Gaussian integral peaked at $x = s - 1$ for $s \to \infty$. The Gaussian approximation for the integrand means the quadratic approximation for its logarithm, and the latter is plotted thin in Figure (187).

This connects all the difference appearances of the number $\pi$ in this section.

### E.6.

Therefore, we have the following great mnemonic formula

$$\mathsf{v}_d = \frac{\pi^{d/2}}{(d/2)!}. \tag{198}$$

For odd dimensions, we should define the factorial using the Gamma functions, and concretely

$$(d/2)! = \Gamma(\tfrac{d}{2} + 1) = \tfrac{d}{2} \tfrac{d-2}{2} \tfrac{d-4}{2} \cdots \tfrac{1}{2} \Gamma(\tfrac{1}{2}) = 2^{-(d+1)/2} d!! \sqrt{\pi}, \quad \text{for } d \text{ odd}. \tag{199}$$

Here $d!!$ means the double factorial of an integer $d$, that is, the product of odd (respectively, even) integers in $\{1, \ldots, d\}$.

### E.7.

Note that from the Stirling formula, we have

$$\operatorname{Vol} B(0, r) \sim \frac{1}{\sqrt{\pi d}} \left( \frac{2\pi e\, r^2}{d} \right)^{d/2}$$

which leads to a very remarkable conclusion: the volume of a ball of *arbitrarily large* fixed radius $r$ goes to 0 as $d \to \infty$ superexponentially fast!

### F. More on $E_8$ and regular $m$-gons
### F.1.

In addition to Coxeter elements C, which have order 30, the group $W(E_8)$ generated by reflections in the roots of the $E_8$ lattice has distinguished conjugacy classes of elements of order 24 and 20; see [24]. In this sections, we will denote a representative of these conjugacy classes by $C_{30}$, $C_{24}$, and $C_{20}$. Their eigenvalues of each $C_m$ are the primitive $m$th roots of unity. There are exactly 8 of those in each case.

Projecting the roots on any of the eigenspaces, one gets the following patterns:

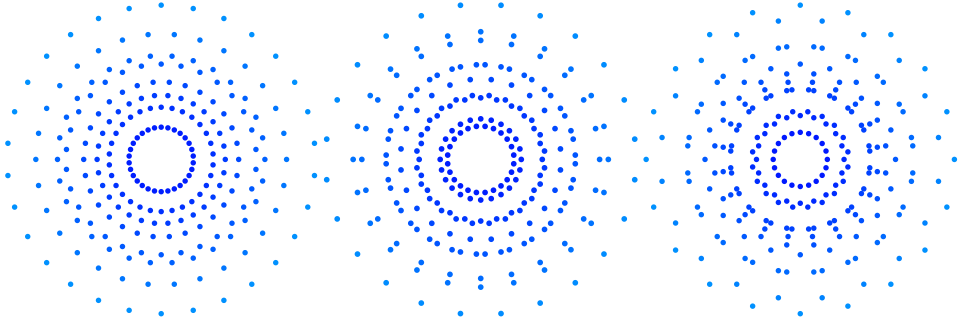

### F.2.

From a slightly different angle, the relation between $E_8$ and regular polygons may be seen as follows. For any $m$, the $m$th roots of unity $\{1, \zeta, \zeta^2, \ldots, \zeta^{m-1}\}$, where $\zeta = \exp\left(\frac{2\pi i}{m}\right)$, are the image of the group $G = \mathbb{Z}/m$ in a 1-dimensional representation. For instance, for $m = 6$ these are the vertices of a regular hexagon:

$$
\begin{array}{c}
\zeta^2 \quad \zeta^1 \\
\zeta^3 \qquad \zeta^0 \\
\zeta^4 \quad \zeta^5
\end{array}
\tag{200}
$$

We can also consider the image of the group ring $\mathbb{Z}G$, that is, the subring

$$\Lambda = \mathbb{Z}[\zeta] \subset \mathbb{C}.$$

This is the set of points that can be obtained by adding and subtracting the vertices of a regular $m$-gon. For instance for $m = 4, 6$ this will be the square lattice $A_1 \oplus A_1$ and the hexagonal lattice $A_2$, respectively.

The powers $1, \zeta, \zeta^2, \ldots$ are linearly independent over $\mathbb{Q}$ until we get to $\zeta^{\phi(m)}$, where $\phi(m)$ is number of residues modulo $m$ that are coprime to $m$, also known as Euler's totient. The number $\zeta^{\phi(m)}$ is an integral linear combination of the numbers $1, \zeta, \ldots, \zeta^{\phi(m)-1}$, given by the coefficients of the cyclotomic polynomial

$$\Psi_m(x) = \prod_{\gcd(i,m)=1}(x - \zeta^i) = x^{\phi(m)} + \cdots \in \mathbb{Z}[x]. \tag{201}$$

See (51) for the explicit form of $\Psi_{30}$. Thus,

$$\Lambda \cong \mathbb{Z}^{\phi(m)} \subset \mathbb{C} \tag{202}$$

as a group under addition. Since we want to construct $E_8$, we focus our attention on the case

$$\phi(m) = 8 \quad \Rightarrow \quad m \in \{15, 20, 24, 30\}. \tag{203}$$

The number 15 here corresponds to the element $\mathsf{C}_{30}^2$. We skip it, since the 15-gon and the 30-gon generate the same $\Lambda$.

Lest the reader imagine $\Lambda \subset \mathbb{C}$ as a lattice, we plot the points $\sum_{i=0}^{29} c_i \zeta^i$, where $\zeta = \exp\left(\frac{\pi i}{15}\right)$ and $c_i \in \{0, 1, 2, 3\}$. $\Lambda$ is a free abelian subgroup of $\mathbb{C}$, but it is not a lattice.

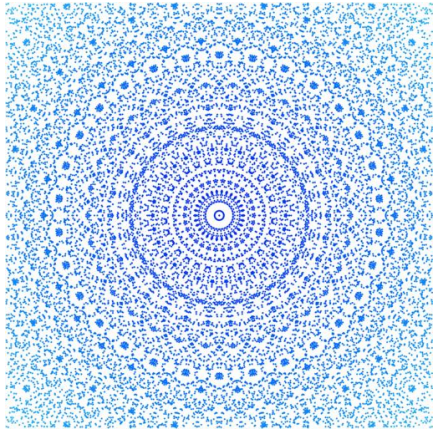

(204)

### F.3.

Let $v \in E_8$ be a vector such that the vectors $\mathsf{C}_m^i v$ span the lattice $E_8$, and consider the diagonal matrix element

$$\eta_v(\zeta^i) = (\mathsf{C}_m^i v, v)_{E_8}. \tag{205}$$

Since $\Psi_m(\mathsf{C}_m) = 0$, this gives a well-defined linear function on $\Lambda$. Put slightly differently, since the eigenvalues of $\mathsf{C}_m$ are the primitive roots of unity, only those Fourier coefficients of $\eta_v$, viewed as a function on $\mathsf{G}$, do not vanish. This makes it well-defined as a function on $\Lambda$. Furthermore, the function (58) being positive-definite, these Fourier coefficients are positive.

For example, let us take the particular Coxeter element constructed in (50) and the following vectors

$$v = \alpha_5, \quad v' = 2e_8, \tag{206}$$

that is, the triple node the Dynkin diagram and twice the last coordinate vector. From the explicit expression

$$\mathsf{C}_{30} = \frac{1}{4}\begin{bmatrix} -1 & -1 & 3 & -1 & -1 & 1 & 1 & -1 \\ 3 & -1 & -1 & -1 & -1 & 1 & 1 & -1 \\ -1 & -1 & -1 & -1 & 3 & 1 & 1 & -1 \\ -1 & 3 & -1 & -1 & -1 & 1 & 1 & -1 \\ -1 & -1 & -1 & -1 & -1 & -3 & 1 & -1 \\ -1 & -1 & -1 & 3 & -1 & 1 & 1 & -1 \\ -1 & -1 & -1 & -1 & -1 & 1 & -3 & -1 \\ -1 & -1 & -1 & -1 & -1 & 1 & 1 & 3 \end{bmatrix}, \tag{207}$$

one can check that the corresponding functions $\eta$ and $\eta'$ are given by

$$\eta(\zeta^i) = \left\lfloor\!\!\left\lfloor 2\cos\left(\frac{\pi i}{15}\right)\right\rfloor\!\!\right\rfloor, \quad \eta'(\zeta^i) = \left\lfloor\!\!\left\lfloor 4\cos\left(\frac{\pi i}{15}\right)\right\rfloor\!\!\right\rfloor, \tag{208}$$

where $\lfloor\!\lfloor x \rfloor\!\rfloor$ denotes the integer between 0 and $x$ that is closest to $x$. These formulas may be illustrated as follows

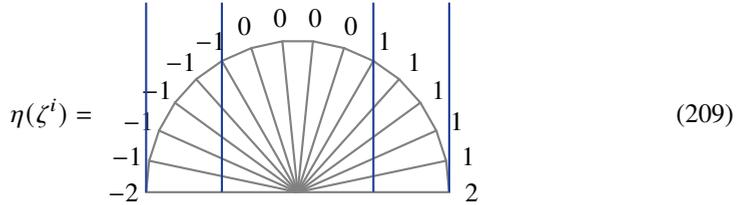

$$\eta(\zeta^i) = \tag{209}$$

and

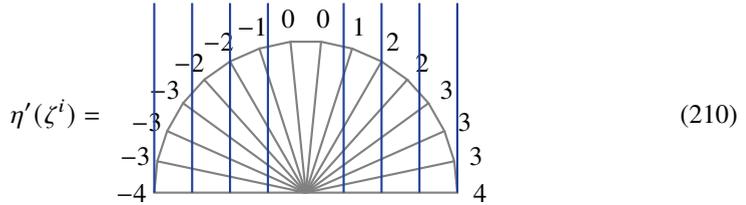

$$\eta'(\zeta^i) = \tag{210}$$

The proximity of these functions to the cosine function can be interpreted as the proximity of the vectors $v$ and $v'$ to the plane in which $\mathsf{C}_{30}$ acts as a rotation by $\pi/15$.

### F.4.

In the style of Appendix B.6, one can turn this construction around as follows. One can check directly that the Fourier coefficients satisfy

$$\widehat{\eta}(j) \text{ is } \begin{cases} > 0, & \gcd(j, m) = 1, \\ = 0, & \gcd(j, m) > 1, \end{cases} \tag{211}$$

and similarly for $\eta'$. We can then *define* the $E_8$ lattice as the group $\Lambda$ with the quadratic form

$$(v, v')_{E_8} = \eta(v\overline{v'}), \quad v, v' \in \Lambda. \tag{212}$$

In fact, the functions

$$\eta_m(\zeta^i) = \left\lfloor\!\!\left\lfloor 2\cos\left(\frac{2\pi i}{m}\right)\right\rfloor\!\!\right\rfloor, \quad \eta'_m(\zeta^i) = \left\lfloor\!\!\left\lfloor 4\cos\left(\frac{2\pi i}{m}\right)\right\rfloor\!\!\right\rfloor, \quad m \in \{20, 24, 30\}, \tag{213}$$

all work and exhibit the $E_8$ lattice as a lattice with an isometry $\mathsf{C}_m$ of the corresponding order. In this realization, the isometry is given by multiplication by $\zeta$. For specific $m$, the numbers 2 and 4 in (213) can be replaced by other even integers.

**Andrei Okounkov**

Andrei Okounkov, Department of Mathematics, University of California, Berkeley, 970 Evans Hall Berkeley, CA 94720–3840, okounkov@math.columbia.edu